\documentclass[10pt,a4paper]{article}

\usepackage{amssymb}
\usepackage{amsmath}
\usepackage{algorithm}     %serve per scrivere lo pseudo codice
\usepackage{algpseudocode} %serve per scrivere lo pseudo codice
\usepackage{subfig}
\usepackage{csquotes}
\usepackage{tikz}
\usetikzlibrary{arrows}
\usepackage{mathtools}
\usepackage{bbm}
\usepackage{amsthm}
\usepackage{pgfplots}
\usepackage{float}
\usepackage{hyperref}
\usepackage{mathbbol}

\usepgfplotslibrary{external}
\newcommand{\vertiii}[1]{{\left\vert\kern-0.25ex\left\vert\kern-0.25ex\left\vert #1 
   \right\vert\kern-0.25ex\right\vert\kern-0.25ex\right\vert}}

\usepackage{geometry}
\makeatletter
\def\blfootnote{\xdef\@thefnmark{$\star$}\@footnotetext}
\makeatother
\newenvironment{Authors}%
  {\begin{center}\begin{bfseries}}%
  {\end{bfseries}\end{center}}
\newenvironment{Addresses}%
  {\begin{flushleft}\begin{itshape}}%
  {\end{itshape}\end{flushleft}}
  \newcommand{\email}[1]{\hspace*{\stretch{1}}\emph{\texttt{#1}}}

 \usepackage{fancyhdr}  
  \fancypagestyle{plain}{
\fancyhead{}
\fancyhead[C]{\hfill Submitted to Springer, January 2024}
 }
\begin{document}

\thispagestyle{plain}

\title{Accelerated construction  of projection-based reduced-order models via incremental approaches}
 \date{}
 \maketitle

 \maketitle
\vspace{-50pt} 
 
\begin{Authors}
Eki Agouzal$^{1,2}$, Tommaso Taddei$^{2}$
\end{Authors}

\begin{Addresses}
$^1$
EDF Lab Paris-Saclay,
7 Boulevard Gaspard Monge, 91120, 
Palaiseau, France.\\
\email{eki.agouzal@edf.fr} \\
$^2$
Univ. Bordeaux, CNRS, Bordeaux INP, IMB, UMR 5251, F-33400 Talence, France. \\ Inria Bordeaux Sud-Ouest, Team MEMPHIS, 33400 Talence, France.
\email{eki.agouzal@inria.fr,tommaso.taddei@inria.fr} \\
\end{Addresses}

\begin{abstract}
We present an accelerated greedy strategy for training of projection-based reduced-order models for parametric steady and unsteady partial differential equations.
Our approach exploits hierarchical approximate proper orthogonal decomposition to speed up the construction of the empirical test space for least-square Petrov-Galerkin formulations, a progressive construction of the empirical quadrature rule based on a warm start of the non-negative least-square algorithm, and a two-fidelity sampling strategy to reduce the number of expensive greedy iterations. We illustrate the performance of our method for two test cases: a two-dimensional compressible inviscid flow past a LS89 blade at moderate Mach number, and a three-dimensional nonlinear mechanics problem to predict the long-time structural response of the standard section of a nuclear containment building  under external loading.
\end{abstract}

\medskip
 
   \section{Introduction}
\label{sec:intro}

In the past few decades, several studies have shown the potential of model order reduction (MOR) techniques to speed up the solution to  many-query and real-time problems, and ultimately enable the use of physics-based three-dimensional  models for design and   optimization, uncertainty quantification, real-time control and monitoring tasks.
The distinctive feature of MOR methods
is the \emph{offline/online computational decomposition}: during the  \emph{offline stage},  high-fidelity (HF)  simulations are employed to  generate an empirical reduced-order approximation of the solution field 
and a parametric reduced-order model (ROM);
during the  \emph{online stage},
the ROM is solved to estimate the solution field and relevant quantities of interest for several parameter values.
Projection-based ROMs (PROMs) rely on the projection of the equations onto a suitable low-dimensional test space.
Successful MOR techniques should hence achieve significant online speedups at  acceptable offline training costs.
This work addresses the reduction of offline training costs  of   PROMs
for  parametric steady and unsteady partial differential equations (PDEs).

We denote by $\mu$ the vector of $p$ model parameters in the 
compact parameter region
$\mathcal{P}\subset \mathbb{R}^p$; given the domain $\Omega\subset \mathbb{R}^d$ ($d=2$ or $d=3$), we introduce the Hilbert spaces $(\mathcal{X}, \|\cdot \|)$ and
$(\mathcal{Y}, \vertiii{\cdot})$ defined over $\Omega$. Then, we consider problems of the form
\begin{equation}
\label{eq:abstractPDE}
{\rm find} \; u_{\mu}\in \mathcal{X} \; : \; 
\mathfrak{R}_{\mu}(u_{\mu}, v) = 0, \quad
\forall \; v\in \mathcal{Y}, 
\;\;\mu\in \mathcal{P},
\end{equation}
where $\mathfrak{R}:\mathcal{X} \times
\mathcal{Y} \times \mathcal{P} \to \mathbb{R}$ is the parametric residual associated with the PDE of interest. We here focus on linear approximations, that is we consider reduced-order approximations of the form
\begin{equation}
\label{eq:linear_ansatz}
\widehat{u}_{\mu} = Z \widehat{\alpha}_{\mu}, 
\end{equation}
where $Z:\mathbb{R}^n \to \mathcal{X}$ is a linear operator, and $n$ is much smaller than the size of the HF model;
$\widehat{\alpha}: \mathcal{P} \to \mathbb{R}^n$ is the vector of generalized coordinates.
As discussed in section \ref{sec:steady_problems}, we here exploit the least-square Petrov-Galerkin (LSPG, 
\cite{carlberg2013gnat}) ROM formulation proposed in \cite{taddei2021space}: the approach relies on the definition of a low-dimensional empirical test space 
$\widehat{\mathcal{Y}} \subset \mathcal{Y}$; furthermore, it relies on the approximation of the HF residual $\mathfrak{R}$ through hyper-reduction 
\cite{ryckelynck2009hyper} to enable fast online calculations of the solution to the ROM. In more detail, we consider an empirical quadrature (EQ) procedure 
\cite{farhat2015structure,yano2019lp} for hyper-reduction: EQ methods recast the problem of finding a sparse quadrature rule to approximate 
$\mathfrak{R}$ as a sparse representation problem and then resort to optimization algorithms to find an approximate solution.
Following \cite{farhat2015structure}, we resort to the non-negative least-square (NNLS) method to find the quadrature rule.

In section \ref{sec:unsteady_problems}, we also consider the application to unsteady problems: to ease the presentation, we consider one-step time discretizations based on the time grid $\{ t^{(k)}  \}_{k=0}^K$. Given $\mu\in \mathcal{P}$, we seek the sequence $\mathbb{u}_{\mu} = 
\{  {u}_{\mu}^{(k)}  \}_{k=0}^K$ such that
\begin{equation}
\label{eq:abstract_unsteadyPDE}
\left\{
\begin{array}{ll}
\displaystyle{
\mathfrak{R}_{\mu}^{(k)}(u_{\mu}^{(k)},
u_{\mu}^{(k-1)}, 
 v) = 0
} &
\forall \, v\in \mathcal{Y} ;\\[3mm]
u_{\mu}^{(0)} = \bar{u}_{\mu}^{(0)};&
 \\[3mm]
\end{array}
\right.
\end{equation}
for all $\mu\in \mathcal{P}$.
As for the steady-state case, we consider linear ansatzs of the form 
$\widehat{u}_{\mu}^{(k)} = Z \widehat{\alpha}_{\mu}^{(k)}$, 
where $Z:\mathbb{R}^n \to \mathcal{X}$ is a linear time-
and parameter-independent operator, and 
$\widehat{\alpha}^{(0)}, \ldots, 
\widehat{\alpha}^{(K)}
: \mathcal{P} \to \mathbb{R}^n$ are obtained by projecting the equations \eqref{eq:abstract_unsteadyPDE} onto a low-dimensional test space.
 To speed up online costs, we also replace the HF residual  in \eqref{eq:abstract_unsteadyPDE} with a rapidly-computable surrogate through the same  hyper-reduction technique considered for steady-state problems.

Following the seminal work by Veroy \cite{veroy2003posteriori}, numerous authors have
resorted to greedy methods to adaptively sample the parameter space, with the ultimate aim of reducing the number of HF simulations performed during the training phase. 
Algorithm \ref{alg:weak_greedy} summarizes the general methodology.
First, we initialize the reduced-order basis (ROB) $Z$ and the ROM based on \emph{a priori} sampling of the parameter space; second, we repeatedly solve the ROM and we estimate the error over a range of parameters $\mathcal{P}_{\rm train}\subset \mathcal{P}$; third, we compute the HF solution for the parameter that maximizes the error indicator; fourth, if the error is above a certain threshold, we update the ROB $Z$ and the ROM, and we iterate; otherwise, we terminate. 
Note that the algorithm depends on an \emph{a posteriori} error indicator $\Delta$: if $\Delta$ is a rigorous  \emph{a posteriori} error estimator, we might apply the termination criterion  directly to the error indicator (and hence save one HF solve).
The   methodology has been extended 
to unsteady problems in \cite{haasdonk2008reduced}:
the method in \cite{haasdonk2008reduced} combines a greedy search driven by an \emph{a posteriori} error indicator with 
proper orthogonal decomposition (POD,
\cite{sirovich1987turbulence,volkwein2011model}) to compress the temporal trajectory.

\begin{algorithm}[H]                      
\caption{: abstract weak greedy algorithm. }     
\label{alg:weak_greedy}     

 \normalsize 

\begin{algorithmic}[1]
\State
Choose 
$\mathcal{P}_{\star} =  \{ \mu^{\star, i}  \}_{i=1}^{n_0}$ and compute the HF solutions
$\mathcal{S}_{\star} =  \{ {u}_{\mu}^{\rm hf}  : \mu \in  \mathcal{P}_{\star}   \} $.
\vspace{3pt}

\For {$n=n_0+1,\ldots,n_{\rm max}$}

\State
Update the ROB ${Z}$ and the ROM.
\vspace{3pt}

\State
Compute the estimate   $\widehat{u}_{\mu}^{\rm hf}$   and evaluate the error indicator 
$\Delta_{\mu}$
for all $\mu \in \mathcal{P}_{\rm train}$.

\vspace{3pt}
 
\State
Compute $u_{\mu^{\star,n}}^{\rm hf}$ for 
 $ \mu^{\star,n} = {\rm arg} \max_{ \mu\in \mathcal{P}_{\rm train}  } \Delta_{\mu}$;
update $\mathcal{P}_{\star}$ and 
$\mathcal{S}_{\star}$. 
\vspace{3pt}

\If{$\| {u}_{\mu^{\star,n}}^{\rm hf} - \widehat{u}_{\mu^{\star,n}}   \| < {\texttt{tol}} 
\| {u}_{\mu^{\star,n}}^{\rm hf}  \|$}

\State
Update the ROB ${Z}$ and the ROM.

\State
\texttt{break}
\EndIf

\EndFor
\end{algorithmic}
\bigskip

\end{algorithm}
 
As observed by several authors, greedy methods enable effective sampling of the parameter space \cite{cohen2015approximation}; however, they suffer from several limitations that might ultimately limit their effectiveness compared to standard \emph{a priori} sampling.
First, Algorithm \ref{alg:weak_greedy} is inherently sequential and cannot hence benefit from parallel architectures.
Second, Algorithm \ref{alg:weak_greedy} requires the solution to the ROM and the evaluation of the error indicator for several parameter values at each iteration of the offline stage; similarly, it requires the update of the ROM --- i.e., the trial and test bases, the reduced quadrature, and possibly the data structures employed to evaluate the error indicator. These observations motivate the development of more effective training strategies for MOR.

We propose an acceleration strategy for Algorithm \ref{alg:weak_greedy} based on three ingredients: 
(i)  a hierarchical construction of the empirical test space for LSPG ROMs based on the hierarchical approximate  POD (HAPOD, 
\cite{himpe2018hierarchical}); (ii) a progressive construction of the empirical quadrature rule based on a warm start of the NNLS algorithm;
(iii) a two-fidelity sampling strategy that is based on the application  of the strong-greedy algorithm 
(see, e.g., \cite[section 7.3]{quarteroni2015reduced})
to a dataset of coarse simulations.
Re (iii), sampling based on coarse simulations is employed to initialize Algorithm \ref{alg:weak_greedy} (cf. Line 1) and ultimately reduce the number of expensive greedy iterations.
We illustrate the performance of our method for two test cases: a two-dimensional  compressible inviscid flow past a LS89 blade at moderate Mach number, and a three-dimensional nonlinear mechanics problem to predict the long-time structural response of 
the standard section of a nuclear containment building (NCB) under external loading.

Our method shares several features with previously-developed techniques.  Incremental POD techniques have been extensively applied to avoid the storage of the full snapshot set for unsteady simulations (see, e.g., \cite{himpe2018hierarchical} and the references therein): here, we adapt the incremental approach described in  \cite{himpe2018hierarchical} to the construction of the test space for LSPG formulations. Chapman  
\cite{chapman2017accelerated}  extensively discussed the parallelization of the NNLS algorithm for MOR applications: we envision that our method can be easily combined with the method in \cite{chapman2017accelerated} to further reduce the training costs. Several authors have also devised strategies to speed up the greedy search through the vehicle of a surrogate error model 
\cite{feng2023multi,paul2015adaptive}; our multi-fidelity strategy extends the work by Barral \cite{barral2023registration} to  unsteady problems and to a more challenging --- from the perspective of the HF solver ---  compressible flow test.
We remark that our approach is similar in scope to the work by Benaceur \cite{benaceur2018progressive} 
that devised a progressive    empirical interpolation method (EIM, 
\cite{barrault2004empirical}) for  hyper-reduction of nonlinear problems.
Finally, we observe that multi-fidelity techniques have been extensively considered for non-intrusive MOR (see, e.g.,
\cite{conti2023multi} and the references therein).

The paper is organized as follows.
In section \ref{sec:steady_problems} we discuss the methodology for steady-state problems: we first review the LSPG hyper-reduced ROM formulation and we review the construction of the quadrature rule; then, we present the accelerated strategy.
In section \ref{sec:unsteady_problems} we extend the method to unsteady problems of the form \eqref{eq:abstract_unsteadyPDE}: for simplicity, we only consider the case of Galerkin ROMs.
Section \ref{sec:numerics} illustrates the performance of the method for  two model problems.
Section \ref{sec:conclusions} draws some conclusions and discusses    future developments.

\section{Accelerated construction of least-square Petrov-Galerkin ROMs for steady problems}
\label{sec:steady_problems}
Section \ref{sec:preliminaries} summarizes  notation that is employed throughout the paper and reviews three  algorithms --- POD, strong-greedy, and  the active set method for NNLS problems --- that are used afterwards.
Section \ref{sec:lspg} reviews the LSPG formulation employed  in the present work and illustrates the strategies employed for the definition of the quadrature rule and the empirical test space.
Section \ref{sec:progressive_EQETS} discusses the incremental strategies proposed in this work, to reduce the costs associated with the construction of the ROM.
Finally, section \ref{sec:multifidelity_sampling} summarizes the multi-fidelity sampling strategy, which is designed to reduce the total number of greedy iterations required by Algorithm \ref{alg:weak_greedy}.

\subsection{Preliminary definitions and tools}
\label{sec:preliminaries}

We denote by 
$\mathcal{T}_{\rm hf} = \left(   \{   {x}_j^{\rm hf}\}_{j=1}^{N_{\rm nd}}, \texttt{T}   \right)$  the HF mesh of the domain $\Omega$ 
with nodes  $\{   {x}_j^{\rm hf}\}_j$ and  connectivity matrix $\texttt{T}$; 
we introduce the elements
$\{ \texttt{D}_k  \}_{k=1}^{N_{\rm e}}$ and the facets
$\{ \texttt{F}_j  \}_{j=1}^{N_{\rm f}}$ of the mesh and we define the open set $\widetilde{\texttt{F}}_j$ in $\Omega$ as the union of the elements of the mesh that share the facet  
$\texttt{F}_j$, with $j=1,\ldots,N_{\rm f}$.
We further denote by $\mathcal{X}_{\rm hf}$ the HF space associated with the mesh $\mathcal{T}_{\rm hf}$ and we define
$N_{\rm hf}:={\rm dim} ( \mathcal{X}_{\rm hf}   )$.

Exploiting the previous definitions, we can introduce the HF residual:
\begin{equation}
\label{eq:HF_residual}
\mathfrak{R}_{\mu}^{\rm hf}(w, v)
=
\sum_{k=1}^{N_{\rm e}} 
r_{k,\mu}^{\rm e} \left(  w|_{\texttt{D}_k},  v|_{\texttt{D}_k}  \right) \; + \;
\sum_{j=1}^{N_{\rm f}} r_{j, \mu}^{\rm f} \left(
 w|_{\widetilde{\texttt{F}}_j}   , v|_{\widetilde{\texttt{F}}_j}
 \right),
\end{equation}
for all  $w,v\in  \mathcal{X}_{\rm hf}$.
To shorten notation, in the remainder, we do not explicitly include the restriction operators in the residual.
The global residual can be viewed as the sum of local element-wise residuals and local facet-wise residuals, which can be evaluated at a cost that is independent of the total number of elements and facets, and is based on local information.
The nonlinear infinite-dimensional statement \eqref{eq:abstractPDE} translates into the high-dimensional problem: find $u_{\mu}^{\rm hf}\in  \mathcal{X}_{\rm hf}  $ such that
\begin{equation}
\label{eq:HF_steady}
\mathfrak{R}_{\mu}^{\rm hf} \left(  u_{\mu}^{\rm hf} , v \right) = 0 \quad
\forall \, v\in   \mathcal{X}_{\rm hf}.
\end{equation}
Since \eqref{eq:HF_steady} is nonlinear, the solution requires to solve a nonlinear system of $N_{\rm hf}$ equations with $N_{\rm hf}$ unknowns.
Towards this end, we here resort to the pseudo-transient continuation (PTC) strategy proposed in \cite{yano2011importance}.

We use the method of snapshots (cf. \cite{sirovich1987turbulence}) to compute POD eigenvalues and eigenvectors. 
Given the snapshot set $\{ u^k \}_{k=1}^{n_{\rm train}} \subset \mathcal{X}_{\rm hf}$ and the inner product $(\cdot, \cdot)_{\rm pod}$, we define the Gramian matrix
$\mathbf{C} \in \mathbb{R}^{n_{\rm train} \times  n_{\rm train}}$,
$\mathbf{C}_{k,k'}= (u^k, u^{k'})_{\rm pod}$, 
and we define the POD eigenpairs    
$\{(\lambda_n, \zeta_n)  \}_{n=1}^{n_{\rm train}}$
as
$$
\mathbf{C} \boldsymbol{\zeta}_i = \lambda_i  \, \boldsymbol{\zeta}_i,
\quad
\zeta_i:=   \sum_{k=1}^{n_{\rm train}} \, \left(   \boldsymbol{\zeta}_i  \right)_k \, u^k,
\quad
i=1,\ldots,n_{\rm train},
$$
with $\lambda_1 \geq \lambda_2 \geq \ldots \lambda_{n_{\rm train}} \geq 0$. In our implementation, we orthonormalize the modes, that is
$( \zeta_n, \zeta_n)_{\rm pod}= 1$ for $n=1,\ldots,n_{\rm train}$.
In the remainder we use notation
\begin{equation}
\label{eq:POD}
[  \{ \zeta_i \}_{i=1}^n ,  \{  \lambda_i \}_{i=1}^{n}
  ]  
=
\texttt{POD} 
\left( 
\{ u^k \}_{k=1}^{n_{\rm train}}, 
n ,
\; \;
(\cdot, \cdot)_{\rm pod}
\right)
\end{equation}
to refer to the application of POD to the snapshot set 
$\{ u^k \}_{k=1}^{n_{\rm train}}$.
The number of modes $n$ can be chosen 
adaptively by ensuring that the retained energy content is above a certain threshold (see, e.g.,  \cite[Eq. (6.12)]{quarteroni2015reduced}).

We further recall the strong-greedy algorithm: 
the algorithm takes as input the snapshot set 
 $\{ u^k \}_{k=1}^{n_{\rm train}}$, an integer $n$, an inner product $(\cdot, \cdot)_{\rm sg}$ and the induced norm $\|\cdot \|_{\rm sg} = \sqrt{(\cdot, \cdot)_{\rm sg}}$, and returns a set of $n$ indices 
 $\mathcal{I}_{\rm sg}\subset \{1,\ldots,n_{\rm train}\}$
 \begin{equation}
 \label{eq:strong_greedy}
 \left[   \mathcal{I}_{\rm sg} \right]  
=
\texttt{strong-greedy} 
\left( 
\{ u^k \}_{k=1}^{n_{\rm train}}, 
n ,
\; \;
(\cdot, \cdot)_{\rm sg}
\right).
 \end{equation}
The dimension  $n$ can be chosen adaptively by ensuring  that the projection error is below a given threshold ${\rm tol}_{\rm sg}>0$,
\begin{equation}
\label{eq:strong_greedy_choice}
n = \min \left\{
n'\in \{1,\ldots,n_{\rm train}\} \, : \,
\max_{k\in \{1,\ldots,n_{\rm train}\}}
\frac{\|  \Pi_{  \mathcal{Z}_{n'}^\perp   }    u^k \|_{\rm sg}}{\|  u^k \|_{\rm sg}}
\leq  {\rm tol}_{\rm sg}
\right\},
\end{equation}
where $\mathcal{Z}_{n'}$ denotes the $n'$-dimensional space obtained after $n'$ steps of the greedy procedure in Algorithm \ref{alg:strong_greedy}, 
$\mathcal{Z}_{n'}^\perp$ is the orthogonal complement of the space $\mathcal{Z}_{n'}$
and
$\Pi_{  \mathcal{Z}_{n'}^\perp   }: \mathcal{X} \to \mathcal{Z}_{n'}^\perp$ is the orthogonal projection operator onto 
$\mathcal{Z}_{n'}^\perp$.

\begin{algorithm}[H]                      
\caption{: strong-greedy algorithm. }     
\label{alg:strong_greedy}     

 \small
\begin{flushleft}
\emph{Inputs:}  $\{ u^k \}_{k=1}^{n_{\rm train}}$ snapshot set, 
$n$ size of the desired reduced space,
$(\cdot, \cdot)_{\rm sg}$ inner product.
\smallskip

\emph{Outputs:} 
$\mathcal{I}_{\rm sg}$ indices of selected snapshots.

\end{flushleft}                      

 \normalsize 

\begin{algorithmic}[1]
\State
Choose 
$\mathcal{Z}_0 = \emptyset$, 
$\mathcal{I}_{\rm sg}= \emptyset$,
set $\mathcal{I} _{\rm train}=\{1,\ldots,n_{\rm train}\}$.
\vspace{3pt}

\For {$i= 1, \ldots, n$}

\State
Compute $i^{\star} = {\rm arg} \max_{i\in \mathcal{I} _{\rm train}  }   \|  \Pi_{ \mathcal{Z}^\perp   } u^i \|_{\rm sg}$.
\vspace{3pt}

\State
Update
$ \mathcal{Z}_i =  \mathcal{Z}_{i-1} \cup {\rm span} \{  u^{i^{\star}} \}$ and
$\mathcal{I}_{\rm sg}=  \mathcal{I}_{\rm sg} \cup   \{ i^{\star} \}$.
\vspace{3pt}
 
\EndFor
\end{algorithmic}
\bigskip

\end{algorithm}
 
We conclude this section by reviewing the active set method \cite{lawson1995solving} that is employed to find a sparse solution to the non-negative least square problem:
\begin{equation}
\label{eq:NNLS_statement}
\min_{\boldsymbol{\rho}\in \mathbb{R}^N} \| \mathbf{G} \boldsymbol{\rho} - \mathbf{b}  \|_2 \quad {\rm s.t.} \;\;
\boldsymbol{\rho}  \geq 0,
\end{equation}
for a given matrix $\mathbf{G}\in \mathbb{R}^{M\times N}$ and a vector $\mathbf{b}\in \mathbb{R}^M$. 
Algorithm \ref{alg:nnls} reviews  the computational procedure. In the remainder, we use notation
\begin{equation}
\label{eq:NNLS_algo}
 \left[ \boldsymbol{\rho}   \right]  
=
\texttt{NNLS} 
\left( 
\mathbf{G},
\mathbf{b} ,
\delta, 
P_0
\right),
 \end{equation}
to refer to the application of Algorithm \ref{alg:nnls}.

Note that the method takes as input a set of indices --- which is  initialized with   the empty set in the absence of prior information --- to initialize the process.
Given the matrix $\mathbf{G}= [\mathbf{g}_1,\ldots,\mathbf{g}_N]$,
the vector $\mathbf{x}\in \mathbb{R}^N$,
 and the set of indices $P=\{ p_i \}_{i=1}^m \subset \{1,\ldots,N\}$, 
we use  notation $\mathbf{G}(:,P):=[\mathbf{g}_{p_1},\ldots,\mathbf{g}_{p_m}    ]\in \mathbb{R}^{M\times m}$   and
$\mathbf{x}(P) = {\rm vec} ( (\mathbf{x})_{p_1},\ldots,  (\mathbf{x})_{p_m}  ) \in \mathbb{R}^m$; we denote by 
$\# P $ the cardinality of the discrete set $P$, and we introduce the complement of $P$ in 
$\{1,\ldots,N  \} $
as 
 $P^{\rm c}=\{1,\ldots,N  \} \setminus P$.
Given the vector $ \mathbf{x} \in \mathbb{R}^N$ and the set of indices $\mathcal{I} \subset \{1,\ldots,N\}$, 
  notation 
$\left[\alpha, {i}^\star\right]= \min_{i\in \mathcal{I} }  
\left( \mathbf{x} \right)_i$ signifies that
$\alpha=  \min_{i\in \mathcal{I} }  
\left( \mathbf{x} \right)_i$ and $ {i}^\star \in \mathcal{I}$ realizes the minimum, $\alpha=\left( \mathbf{x} \right)_{ {i}^\star}$.
The constant $\epsilon>0$ is intended to avoid division by zero and is set to $2^{-1022}$.
The  computational cost of Algorithm \ref{alg:nnls} is dominated by the cost to repeatedly solve the least-square problem at Line 11: in section \ref{sec:numerics}, we hence report the total number of least-square solves 
needed to achieve
 convergence (cf. output $it$ in Algorithm \ref{alg:nnls}).

\begin{algorithm}[H]                      
\caption{: active set method for \eqref{eq:NNLS_statement}. }     
\label{alg:nnls}     

 \small
\begin{flushleft}
\emph{Inputs:}  
$\mathbf{G} \in \mathbf{R}^{M\times N}$,
$\mathbf{b} \in \mathbb{R}^M$,
$\delta >0$,
$P_0$.
\smallskip

\emph{Output:} 
$\boldsymbol{\rho}$  approximate solution to \eqref{eq:NNLS_statement},
$it$ number of iterations to meet convergence criterion.

\end{flushleft}                      

\normalsize 

\begin{algorithmic}[1]
\State
Choose 
$\boldsymbol{\rho}= 0$, 
$\mathbf{w} = \mathbf{C}^\top \boldsymbol{\rho}$,
$P=P_0$, $it=0$.
\vspace{3pt}

\While {\texttt{true}}

\State
Compute $\mathbf{r} = \mathbf{G}(:,P) \mathbf{x}(P)  -  \mathbf{b}$.
\vspace{3pt}

\If {$\# P = N$ or $\|  \mathbf{r} \|_2 \leq 
\delta \|  \mathbf{b} \|_2$}
\State
\texttt{break}
\EndIf

\State
Set 
$ i^{\star} = {\rm arg} \max_{j\notin P} (\mathbf{w})_j$.
\medskip

\State
Set 
$P = P \cup \{  i^{\star}  \}$.
\medskip

 \While {\texttt{true}}
 
 \State $it=it+1$
 
 \State Define 
 $\mathbf{z}\in \mathbb{R}^N$ s.t.
   $\mathbf{z}(P^{\rm c}) = 0$, 
      $\mathbf{z}(P) = \mathbf{G}(:,P) ^\dagger \mathbf{d}$
     
\If {$\mathbf{z}\geq 0$} 
 \State Set 
 $\mathbf{x} = \mathbf{z}$,   $\mathbf{w}=\mathbf{G}^\top (  \mathbf{b} -  \mathbf{G} \mathbf{x} )$
 
\State
\texttt{break}
\EndIf
 
 \State
$\mathcal{I} = \{i\in \{1,\ldots,N\}: (\mathbf{z})_i < 0   \}$.

  \State
$\left[\alpha, {i}^\star\right]= \min_{i\in \mathcal{I} }  
\frac{  (\mathbf{x})_i    }{(\mathbf{x})_i  -  (\mathbf{z})_i + \epsilon}$.

\State
$P = P \setminus \{  {i}^\star  \}$.

\State
 $\mathbf{x}  =   \mathbf{x}  - \alpha (\mathbf{x} - \mathbf{z})$.

 \EndWhile
 
\EndWhile
\end{algorithmic}
\bigskip

\end{algorithm}

\subsection{Least-square Petrov-Galerkin formulation}
\label{sec:lspg}

 Given the reduced-order basis (ROB) 
 $Z=[\zeta_1,\ldots,\zeta_n]:\mathbb{R}^n \to \mathcal{X}_{\rm hf}$, following \cite{barral2023registration}, we consider the LSPG formulation  
 \begin{subequations}
\label{eq:MOR_intro}
\begin{equation}
\label{eq:MOR_intro_a}
\widehat{u}_{\mu} = Z  \widehat{{\alpha}}_{\mu}
\quad
{\rm with} \;\;
\widehat{{\alpha}}_{\mu} \in  {\rm arg} \min_{ {\alpha} \in \mathbb{R}^n  }
\max_{\psi \in \widehat{\mathcal{Y}}} \;
\frac{\mathfrak{R}_{\mu}^{\rm eq}(
Z  {\alpha}, \psi)}{\vertiii{\psi}}
\end{equation}
where $ \widehat{\mathcal{Y}} = {\rm span} \{  \psi_i \}_{i=1}^m$ is a suitable test space that is chosen below and the empirical residual  $\mathfrak{R}_{\mu}^{\rm eq}$ satisfies
\begin{equation} 
\label{eq:MOR_intro_b}
\mathfrak{R}_{\mu}^{\rm eq}(
w, v)
=
\sum_{k=1}^{N_{\rm e}} 
\rho_{k}^{\rm eq,e}
\,
r_{k,\mu}^{\rm e}( w,   v  ) \; + \;
\sum_{j=1}^{N_{\rm f}} 
\rho_{k}^{\rm eq,f}
r_{j,\mu}^{\rm f}(
w,   v  ) 
\end{equation}
where
$\boldsymbol{\rho}^{\rm eq,e} \in \mathbb{R}_+^{N_{\rm e}}$ and
$\boldsymbol{\rho}^{\rm eq,f} \in \mathbb{R}_+^{N_{\rm f}}$ 
are sparse vectors of non-negative weights. 

Provided that 
$ \{  \psi_i \}_{i=1}^m$ is an orthonormal basis of 
$\widehat{\mathcal{Y}}$, we can rewrite 
\eqref{eq:MOR_intro_a} as
\end{subequations}
\begin{equation}
\label{eq:LSPG_algebraic}
\min_{  
{\alpha} \in \mathbb{R}^n  }
\left\|
\boldsymbol{\mathfrak{R}}_{\mu}^{\rm eq}(
 {\alpha} )
\right\|_2,
\quad
{\rm with} 
\;\;
\left(
\boldsymbol{\mathfrak{R}}_{\mu}^{\rm eq}(
 {\alpha} )
\right)_i =
\mathfrak{R}_{\mu}^{\rm eq}
(Z  {\alpha}, \psi_i),
\;\; i=1,\ldots,m,
\end{equation}
which can be efficiently solved using the Gauss-Newton method (GNM). Note that \eqref{eq:LSPG_algebraic} does not explicitly  depend on the choice of the test norm $\vertiii{\cdot}$: dependence on the norm is implicit in the choice  of 
 $ \{  \psi_i \}_{i=1}^m$.
 Formulation \eqref{eq:MOR_intro}-\eqref{eq:LSPG_algebraic} depends on the choice of the test space $\widehat{\mathcal{Y}}$ and the empirical weights 
$\boldsymbol{\rho}^{\rm eq,e}, \boldsymbol{\rho}^{\rm eq,f}$: in the remainder of this section, we address the construction of these ingredients. 

We remark that 
Carlberg \cite{carlberg2013gnat} considered  
a different projection method, which is based on
the minimization of  the Euclidean  norm of the discrete residual:
due to the particular choice of the test norm, the approach of \cite{carlberg2013gnat} does not require the explicit construction of the empirical test space.
We further observe that Yano and collaborators
\cite{yano2019discontinuous,du2022efficient} 
have considered different formulations of the empirical residual   \eqref{eq:MOR_intro_b}.
A thorough comparison
between
different 
projection methods and different 
hyper-reduction techniques  is beyond the scope of the present study.
 
\subsubsection{Construction of the empirical test space}
\label{sec:empirical_space}
As discussed  in \cite{taddei2021space}, the test space $\widehat{\mathcal{Y}}$ should approximate the Riesz representers of the functionals associated with the action of the Jacobian on the elements of the trial ROB. Given the snapshot set of HF solutions
$\{ 
u_{\mu}^{\rm hf} : \mu\in \mathcal{P}_{\rm train}\}$ with 
$\mathcal{P}_{\rm train}=\{\mu^k \}_{k=1}^{n} \subset \mathcal{P}$ and the ROB 
$\{ \zeta_i \}_i$, we hence apply POD to the test snapshot set 
$\mathcal{S}_n^{\rm test} :=
\left\{ {\Psi}_{k,i}: k=1,\ldots,n_{\rm train} ,i=1,\ldots,n
\right\}$, where
\begin{equation}
\label{eq:test_space_snapshot}
(({\Psi}_{k,i}, v  ))
\, = \,
\mathfrak{J}_{\mu^k}^{\rm hf}\left[  \widehat{u}_{\mu^k}^{\rm hf}  \right]
(  {\zeta}_i, v ),
\quad
\forall \, v\in \mathcal{X}_{\rm hf},
\end{equation}
where 
$\mathfrak{J}_{\mu}^{\rm hf} [w]: \mathcal{X}_{\rm hf}\times  \mathcal{X}_{\rm hf} \to \mathbb{R}$ denotes 
 the Fr{\'e}chet derivative of the HF residual at $w$.
It is also useful to provide an approximate computation of the right-hand side of \eqref{eq:test_space_snapshot},
\begin{equation}
\label{eq:test_space_snapshot_FD}
(({\Psi}_{k,i}, v  ))
\, \approx \,
\frac{1}{\epsilon}
\left(
\mathfrak{R}_{\mu^k}^{\rm hf} 
(    {u}_{\mu^k}^{\rm hf} + \epsilon \zeta_i, v )
\, - \,
\mathfrak{R}_{\mu^k}^{\rm hf} 
(   {u}_{\mu^k}^{\rm hf} , v )
\right)
\quad
\forall \, v\in \mathcal{X}_{\rm hf},
\end{equation}
with $|\epsilon| \ll 1$.
The evaluation of the right-hand side of 
 \eqref{eq:test_space_snapshot_FD}   involves the   computation of the residual at ${u}_{\mu^k}^{\rm hf} + \epsilon \zeta_i$ for $i=1,\ldots,n$;
 on the other hand,
 the evaluation of 
 \eqref{eq:test_space_snapshot} requires the computation of  the Jacobian matrix at 
 ${u}_{\mu^k}^{\rm hf}$ and the post-multiplication by the algebraic counterpart of $Z$.
Both
 \eqref{eq:test_space_snapshot}  and 
 \eqref{eq:test_space_snapshot_FD}
--- which are equivalent  in the limit $\epsilon\to 0$---
are  used in the incremental approach of section \ref{sec:progressive_EQETS}. 
 
\subsubsection{Construction of the empirical quadrature rule}
\label{sec:empirical_quadrature}
Following \cite{yano2019lp}, we seek 
$\boldsymbol{\rho}^{\rm eq,e} \in \mathbb{R}_+^{N_{\rm e}}$ and
$\boldsymbol{\rho}^{\rm eq,f} \in \mathbb{R}_+^{N_{\rm f}}$ in
\eqref{eq:MOR_intro_b}  
 such that
 \begin{enumerate}
 \item[(i)]
 (\emph{efficiency constraint})
the number of nonzero entries in 
 $\boldsymbol{\rho}^{\rm eq,e},\boldsymbol{\rho}^{\rm eq,f}$,  
 $ \texttt{nnz} (  \boldsymbol{\rho}^{\rm eq,e})$ and
 $\texttt{nnz} (   \boldsymbol{\rho}^{\rm eq,f} )$,   is as small as possible;
 \item[(ii)]
 (\emph{constant function constraint})
  the constant function is approximated correctly in $\Omega$
 \begin{equation}
 \label{eq:constant_function_constraint}
\Big|
\sum_{k=1}^{N_{\rm e}} \rho_k^{\rm eq,e} | \texttt{D}_k  |
\,-\,
| \Omega | 
 \Big|
 \ll 1,
 \quad
\Big|
\sum_{j=1}^{N_{\rm f}} \rho_j^{\rm eq,f} | \texttt{F}_j  |
\,-\,
\sum_{j=1}^{N_{\rm f}}  | \texttt{F}_j  |
 \Big|
 \ll 1; 
\end{equation}
\item[(iii)]
(\emph{manifold accuracy constraint})
for all $\mu \in \mathcal{P}_{\rm train,eq} = \{  \mu^k \}_{k=1}^{n_{\rm train}+n_{\rm train,eq}}$, the  empirical residual satisfies
\begin{subequations} 
\label{eq:accuracy_constraint}
\begin{equation}
\Big \|
   {\boldsymbol{\mathfrak{R}}}_{\mu}^{\rm hf}
(  {\alpha}_{\mu}^{\rm train}   )   
\, - \,
   {\boldsymbol{\mathfrak{R}}}_{\mu}^{\rm eq}
( {\alpha}_{\mu}^{\rm train}   )   
 \Big \|_2
 \ll 1,
\end{equation}
where 
${\boldsymbol{\mathfrak{R}}}_{\mu}^{\rm hf}$  corresponds to substitute
$\rho_1^{\rm eq,e} = \ldots = \rho_{N_{\rm e}}^{\rm eq,e} = 
\rho_1^{\rm eq,f} = \ldots = \rho_{N_{\rm f}}^{\rm eq,f} = 1$ in
\eqref{eq:MOR_intro_b} 
 and ${\alpha}_{\mu}^{\rm train}$ satisfies
\begin{equation}
{\alpha}_{\mu}^{\rm train} = 
\left\{
\begin{array}{ll}
\displaystyle{ {\rm arg} \min_{ {\alpha} \in \mathbb{R}^n} \; 
\|   Z  {\alpha}    - {u}_{\mu}^{\rm hf} \|,}
 & {\rm if} \; \mu \in \mathcal{P}_{\rm train} ; \\[3mm]
\displaystyle{ {\rm arg} \min_{{\alpha} \in \mathbb{R}^n} \; 
\| {\boldsymbol{\mathfrak{R}}}_{\mu}^{\rm hf}
(  {\alpha}   )   
\|_2,}
  & {\rm if} \; \mu \notin \mathcal{P}_{\rm train} ; \\
\end{array}
\right.
\end{equation}
and 
$\mathcal{P}_{\rm train} = \{  \mu^k \}_{k=1}^{n_{\rm train}}$ is the set of parameters for which the HF solution is available.
In the remainder,  we use $\mathcal{P}_{\rm train,eq} = \mathcal{P}_{\rm train}$.
\end{subequations}
 \end{enumerate}
 
 By tedious but straightforward calculations,  we find that 
 \begin{subequations}
 \label{eq:EQ_explained}
 \begin{equation}
\label{eq:EQ_explained_a}
\mathfrak{R}_{\mu^k}^{\rm eq}
(Z  {\alpha}, \psi_i)
=
\mathbf{G}_{k,i}^{\rm e}  \cdot  \boldsymbol{\rho}^{\rm eq,e}
+
\mathbf{G}_{k,i}^{\rm f}    \cdot   \boldsymbol{\rho}^{\rm eq,f}
\end{equation}
for some row vectors $\mathbf{G}_{k,i}^{\rm e} \in \mathbb{R}^{1\times N_{\rm e}}$
and
$\mathbf{G}_{k,i}^{\rm f} \in \mathbb{R}^{1\times N_{\rm f}}$, 
$k=1,\ldots,n_{\rm train}$ and $i=1,\ldots,m$. 
Therefore, if we define
\begin{equation}
\label{eq:EQ_explained_b}
\mathbf{G}^{\rm cnst,e} =
\left[ |\texttt{D}_1|, \ldots,  |\texttt{D}_{N_{\rm e}}| \right],
\quad
\mathbf{G}^{\rm cnst,f} =
\left[ |\texttt{F}_1|, \ldots,  |\texttt{F}_{N_{\rm f}}| \right],
\end{equation}
we   find  that the constraints \eqref{eq:constant_function_constraint} and \eqref{eq:accuracy_constraint} can be expressed in the algebraic form
\begin{equation}
\label{eq:EQ_explained_c}
\| \mathbf{G}  \left(   \boldsymbol{\rho}^{\rm eq} -    \boldsymbol{\rho}^{\rm hf}  \right) \|_2 \ll 1,
\quad
{\rm with} \;\;
\mathbf{G}   = \left[
\begin{array}{ll}
\mathbf{G}_{1,1}^{\rm e}  &  \mathbf{G}_{1,1}^{\rm f} \\ 
\vdots & \vdots \\
\mathbf{G}_{n_{\rm train},m}^{\rm e}  &  \mathbf{G}_{n_{\rm train},m}^{\rm f} \\ 
\mathbf{G}^{\rm cnst,e}  &  0 \\ 
0 &  \mathbf{G}^{\rm cnst,f}  \\ 
\end{array}
\right],
\;\;
\boldsymbol{\rho}^{\rm eq}=
\left[
\begin{array}{c}
\boldsymbol{\rho}^{\rm eq,e} \\
\boldsymbol{\rho}^{\rm eq,f} \\
\end{array}
\right],
\end{equation}
and $\boldsymbol{\rho}^{\rm hf} = [1,\ldots,1]^\top$.
 \end{subequations}
 
 In conclusion, the problem of finding the sparse weights 
$\boldsymbol{\rho}^{\rm eq,e}, \boldsymbol{\rho}^{\rm eq,f}$ can be recast as a sparse representation problem 
\begin{equation}
\label{eq:sparse_representation}
\min_{ \boldsymbol{\rho}\in \mathbb{R}_+^{N_{\rm e}+N_{\rm f}}    }
\texttt{nnz} (\boldsymbol{\rho} ),
\quad
\| \mathbf{G}      \boldsymbol{\rho}^{\rm eq} -    
\mathbf{b}    \|_2\leq \delta \|
\mathbf{b}    \|_2,
\end{equation}
where 
$\texttt{nnz} (\boldsymbol{\rho} )$ is the number of non-zero entries in
the vector $\boldsymbol{\rho}$,  
for some user-defined tolerance $\delta>0$, and
$\mathbf{b}   = \mathbf{G}   \boldsymbol{\rho}^{\rm hf}$.
Following \cite{farhat2015structure}, we resort to the NNLS algorithm discussed in section \ref{sec:preliminaries} (cf. Algorithm \ref{alg:nnls}) to find an approximate solution to \eqref{eq:sparse_representation}.

\subsubsection{A posteriori error indicator}
\label{sec:greedy_search}
The final element of the formulation is the \emph{a posteriori} error indicator that is employed for the parameter exploration.
We here consider  the residual-based error indicator
(cf. \cite{barral2023registration}),
\begin{equation}
\label{eq:error_indicator}
\Delta: \mu\in \mathcal{P} \mapsto
\sup_{v\in \mathcal{X}_{\rm hf}} 
\frac{\mathfrak{R}_{\mu}^{\rm hf} (\widehat{u}_{\mu}, v)}{\vertiii{v}}.
\end{equation}
Note that the evaluation of \eqref{eq:error_indicator} requires the solution to a linear system of size $N_{\rm hf}$: it is hence ill-suited for real-time online computations; nevertheless, in our experience the offline cost associated with the evaluation of \eqref{eq:error_indicator} is 
comparable with the cost that is needed to solve the ROM --- clearly, this is related to the size of the mesh and is hence strongly problem-dependent. 

\subsubsection{Overview of the computational procedure}
\label{sec:overview_standard}

Algorithm \ref{alg:rom_construction_steady}
  provides a detailed summary of the construction of the ROM at each step of Algorithm \ref{alg:weak_greedy} (cf. Line 7).
Some comments are in order.
The cost of Algorithm \ref{alg:rom_construction_steady} is dominated by the assembly of the test snapshot set and by the solution to the NNLS problem. We also notice that the storage of $\mathcal{S}_n^{\rm test}$ scales with $\mathcal{O}(n^2)$ and is hence the dominant memory cost of the offline procedure.

\begin{algorithm}[H]                      
\caption{: generation of the ROM. }     
\label{alg:rom_construction_steady}     

 \small
\begin{flushleft}
\emph{Inputs:}  
$Z=[\zeta_1,\ldots,\zeta_{n-1}]$ actual ROB,
${\rm tol}>0$ tolerance for Algorithm \ref{alg:nnls},
$u_{\mu^n}^{\rm hf}$ new snapshot,
$m$ size of the test space
\smallskip

\emph{Outputs:} 
$Z=[\zeta_1,\ldots,\zeta_{n}]$   new ROB,
ROM for the generalized coordinates.
\end{flushleft}                      

\normalsize 

\begin{algorithmic}[1]
\State
Apply Gram-Schmidt orthogonalization to define the new ROB $Z$.
\vspace{3pt}

\State
Assemble the test snapshot set $\mathcal{S}_n^{\rm test}$ using 
 \eqref{eq:test_space_snapshot}.
\vspace{3pt}

\State
Apply POD to find $\widehat{\mathcal{Y}}$:
$  \{ \psi_i \}_{i=1}^m 
=
\texttt{POD} 
\left( 
\mathcal{S}^{\rm test}, 
m,
\; \;
((\cdot, \cdot))
\right)
$.
\vspace{3pt}

\State
Assemble the matrix and vector $\mathbf{G}, \mathbf{b}$ (cf. \eqref{eq:EQ_explained}).
\vspace{3pt}

\State
Solve the NNLS problem (cf. Algorithm \ref{alg:nnls}).
\vspace{3pt}
 
\State
Update the ROM data structures.
\vspace{3pt}

\end{algorithmic}

\end{algorithm}

\subsection{Progressive construction of empirical test space and quadrature}
\label{sec:progressive_EQETS}

\subsubsection{Empirical test space}
\label{sec:progressive_EQETS_a}
By reviewing Algorithms \ref{alg:weak_greedy} and \ref{alg:rom_construction_steady}, we notice that the test snapshot set $\mathcal{S}_n^{\rm test}$ satisfies
\begin{equation}
\label{eq:test_space_incremental}
\mathcal{S}_n^{\rm test} = \mathcal{S}_{n-1}^{\rm test}
\cup \left\{  \Psi_{k,n} \right\}_{k=1}^{n-1} \cup \left\{  \Psi_{n,i} \right\}_{i=1}^n;
\end{equation}
therefore, at each iteration, it suffices to solve $2n-1$ 
--- as opposed to $n^2$ ---
Riesz problems of the form \eqref{eq:test_space_snapshot} to define $\mathcal{S}_n^{\rm test}$.
As  in \cite{barral2023registration}, we rely on Cholesky factorization
with   fill-in reducing permutations of rows and columns, to reduce the cost of solving the Riesz problems.
In the numerical experiments, we rely on
\eqref{eq:test_space_snapshot},
which involves the assembly of the Jacobian matrix,
 to compute
 $\left\{  \Psi_{n,1} \right\}_{k=1}^n$, 
 while we consider the finite difference  approximation \eqref{eq:test_space_snapshot_FD} to compute
$ \left\{  \Psi_{k,n} \right\}_{k=1}^{n-1}$ with $\epsilon=10^{-6}$.

In order to lower the memory costs associated with the storage of the test snapshot set $\mathcal{S}_n^{\rm test}$ and also the cost of performing POD, we consider a hierarchical approach to construct the test space $\widehat{\mathcal{Y}}$.
In this work, we apply the (distributed) hierarchical approximate proper orthogonal decomposition: HAPOD is related to incremental singular value decomposition \cite{brand2003fast} and guarantees near-optimal performance with respect to the standard POD (cf. \cite{himpe2018hierarchical}).
Given the POD space $\{ \psi_i \}_{i=1}^m$ such that $(( \psi_i , \psi_j )) = \delta_{i,j}$ and the corresponding eigenvalues $\{ \lambda_i \}_{i=1}^m$,  and the new set  of snapshots 
$\left\{  \Psi_{k,n} \right\}_{k=1}^{n-1} \cup \left\{  \Psi_{n,i} \right\}_{i=1}^n$, HAPOD corresponds to applying POD to a suitable snapshot set that combines information from current and previous iterations,
\begin{equation}
\label{eq:HAPOD}
[  \{ \psi_i^{\rm new} \}_{i=1}^{m_{\rm new}} ,  \{  \lambda_i^{\rm new} \}_{i=1}^{m_{\rm new}}
  ]  
=
\texttt{POD} 
\left( 
\mathcal{S}_n^{\rm incr},
\; \;
m_{\rm new} ,
\; \;
((\cdot, \cdot))
\right),
\end{equation}
with 
$\mathcal{S}_n^{\rm incr}:=
\{ \sqrt{\lambda_i}    \psi_i  \}_{i=1}^{m}  \cup
\left\{  \Psi_{k,n} \right\}_{k=1}^{n-1} \cup \left\{  \Psi_{n,i} \right\}_{i=1}^n.$
Note that the modes $\{    \psi_i  \}_{i=1}^{m} $ are scaled to properly take into account the energy content of the modes in the subsequent iterations.
Note also that the storage cost of the method scales with  $m + 2n-1$, which is much lower than $n^2$, provided that $m\ll n^2$. 
Note also  that the POD spaces $\{ \widehat{\mathcal{Y}}_n \}_n$, which are generated at each iteration of   Algorithm \ref{alg:weak_greedy}
 using \eqref{eq:HAPOD}, are not nested.

\subsubsection{Empirical quadrature}
\label{sec:progressive_EQETS_b}

Exploiting \eqref{eq:EQ_explained}, it is straightforward to verify that if the test spaces are nested --- that is, $\widehat{\mathcal{Y}}_{n-1} \subset \widehat{\mathcal{Y}}_n$, then the EQ matrix   $\mathbf{G}_n$ 
at iteration $n$
satisfies
\begin{equation}
\label{eq:incrementalEQ_mat}
\mathbf{G}_n = \left[
\begin{array}{l}
\mathbf{G}_{n-1} \\
\mathbf{G}_n^{\rm new}
\end{array}
\right]
\end{equation}
where $\mathbf{G}_n^{\rm new}$ has 
$k = n \times {\rm dim}( \widehat{\mathcal{Y}}_n \setminus \widehat{\mathcal{Y}}_{n-1}  )
+ {\rm dim}( \widehat{\mathcal{Y}}_{n-1} )$ rows.
We hence observe that we can reduce the cost of assembling the EQ matrix 
$\mathbf{G}_n$ by exploiting \eqref{eq:incrementalEQ_mat}, provided that we rely on nested test spaces: since, in our experience, the cost of assembling $\mathbf{G}_n $ is negligible, the use of non-nested test spaces does not hinder offline performance.

On the other hand, due to the strong link between consecutive EQ problems that are solved at each iteration of the greedy procedure, we propose to initialize the NNLS algorithm \ref{alg:nnls} using the solution from the previous time step, that is 
\begin{equation}
\label{eq:initializationEQP}
P = \left\{ i\in \{ 1,\ldots,N_{\rm e}+ N_{\rm f}\} : (\boldsymbol{\rho}^{{\rm eq}, (n-1)})_i   \neq  0\right\}.
\end{equation}
We  provide extensive numerical investigations to assess the effectiveness of this choice.

\subsubsection{Summary of the incremental weak-greedy algorithm}

Algorithm \ref{alg:rom_construction_steady_incr} reviews  the full  incremental generation of the ROM. In the numerical experiments, we set $m=2n$: this implies that the storage of $\mathcal{S}^{\rm incr}$ scales with  
$4 n - 3$ as opposed to $n^2$.

%
%\ref{alg:rom_construction_steady}  
%
%
%\ref{alg:rom_construction_steady_incr}  

\begin{algorithm}[H]                      
\caption{: incremental generation of the ROM. }     
\label{alg:rom_construction_steady_incr}     

 \small
\begin{flushleft}
\emph{Inputs:}  
$Z=[\zeta_1,\ldots,\zeta_{n-1}]$ actual ROB,
$\delta >0$ tolerance for Algorithm \ref{alg:nnls},
$u_{\mu^n}^{\rm hf}$ new snapshot,
$\{ (\psi_i^{{\rm old}}  \, \lambda_i^{\rm old}  )  \}_{i=1}^{m'}$ POD eigenpairs for test space,
$m$ size of the test space,
$P^{(n-1)}$ initial condition for Algorithm \ref{alg:nnls}.
\smallskip

\emph{Outputs:} 
$Z=[\zeta_1,\ldots,\zeta_{n}]$   new ROB,
ROM for the generalized coordinates,
$\{ (\psi_i ,  \lambda_i)   \}_{i=1}^m$ 
new POD eigenpairs for test space,
$P^{(n)}$ initial condition for Algorithm \ref{alg:nnls} at subsequent iteration.
\end{flushleft}                      
\medskip

\normalsize 

\begin{algorithmic}[1]
\State
Apply Gram-Schmidt orthogonalization to define the new ROB $Z$.
\vspace{3pt}

\State
Assemble the test snapshot set $\mathcal{S}_n^{\rm incr}$ in \eqref{eq:HAPOD}.
\vspace{3pt}

\State
Apply POD to find $\widehat{\mathcal{Y}}$, 
$[  \{ \psi_i  \}_{i=1}^{m} ,  \{  \lambda_i  \}_{i=1}^{m}
  ]  
=
\texttt{POD} 
\left( 
\mathcal{S}_n^{\rm incr},
\; \;
m ,
\; \;
((\cdot, \cdot))
\right)$.
\vspace{3pt}

\State
Assemble the matrix and vector $\mathbf{G}, \mathbf{b}$ (cf. \eqref{eq:EQ_explained}).
\vspace{3pt}

\State
Solve the NNLS problem (cf. Algorithm \ref{alg:nnls}) with initial condition given by $P^{(n-1)}$.
\vspace{3pt}
 
\State
Update the ROM data structures and compute  $P$ using \eqref{eq:initializationEQP}.
\vspace{3pt}

\end{algorithmic}
\end{algorithm}

\subsection{Multi-fidelity sampling}
\label{sec:multifidelity_sampling}

The incremental strategies of section \ref{sec:progressive_EQETS} do not affect the total number of greedy iterations required by Algorithm \ref{alg:weak_greedy} to achieve the desired accuracy; furthermore, they do not address the reduction of the cost of the HF solves. Following \cite{barral2023registration}, we here propose to resort to coarser simulations to learn the initial training set $\mathcal{P}_{\star}$ (cf. Line 1 Algorithm \ref{alg:weak_greedy}) and to initialize the HF solver. Algorithm  \ref{alg:multifidelity_greedy} summarizes the computational procedure.

\begin{algorithm}[H]                      
\caption{: multi-fidelity weak-greedy algorithm. }     
\label{alg:multifidelity_greedy}     

\begin{algorithmic}[1]
\State
Define two meshes   $\mathcal{T}_{\rm hf}^0$ and $\mathcal{T}_{\rm hf}$, where 
$\mathcal{T}_{\rm hf}^0$ is coarser than $\mathcal{T}_{\rm hf}$.
\medskip

\State
Generate a ROM for the HF model associated with the coarse mesh
$\mathcal{T}_{\rm hf}^0$ (e.g. using Algorithm \ref{alg:weak_greedy}),
 $Z^0, \widehat{\alpha}^0: \mathcal{P} \to  \mathbb{R}^n$.
\medskip

\State
Estimate the solution for all $\mu\in  \mathcal{P}_{\rm train}$ and store the generalized coordinates
$\{ \widehat{\alpha}_{\mu}^0  :  \mu\in  \mathcal{P}_{\rm train} \}$.
\medskip

\State
Apply Algorithm \ref{alg:strong_greedy} to 
$\{ \widehat{\alpha}_{\mu}^0  :  \mu\in  \mathcal{P}_{\rm train} \}$
to obtain the initial sample $\mathcal{P}_\star$.
\medskip

\State
Apply Algorithm \ref{alg:weak_greedy} to the HF model associated with the fine mesh $\mathcal{T}_{\rm hf}$; use $\mathcal{P}_\star$ to initialize the greedy method, and use the estimate $\mu \mapsto Z^0 \widehat{\alpha}_\mu^0$  to initialize the HF solver. 

\end{algorithmic}
\end{algorithm}

Some comments are in order.
\begin{itemize}
\item
In the numerical experiments, we choose the cardinality $n_0$ of $\mathcal{P}_\star$ according to \eqref{eq:strong_greedy_choice}  with ${\rm tol}_{\rm sg} =\texttt{tol}$ (cf. Algorithm \ref{alg:weak_greedy}). Note that, since the strong greedy algorithm is applied to the generalized coordinates of the coarse ROM, $n_0$ cannot exceed the size $n$ of the ROM.
\item
We observe that increasing the cardinality of $\mathcal{P}_\star$  ultimately leads to a reduction of the number of sequential greedy iterations; it hence enables much more effective parallelization of the offline stage.
Note also that Algorithm \ref{alg:multifidelity_greedy} can be coupled with parametric mesh adaptation tools to build an effective problem-aware mesh $\mathcal{T}_{\rm hf}$ (cf. \cite{barral2023registration}); in this work, we do not exploit this feature of the method.
\item
We observe that the multi-fidelity algorithm \ref{alg:multifidelity_greedy} critically depends on the choice of the coarse  grid $\mathcal{T}_{\rm hf}^0$: an excessively coarse mesh $\mathcal{T}_{\rm hf}^0$ might undermine the quality of the initial sample $\mathcal{P}_\star$ and of the initial condition for the HF solve, while an excessively fine mesh  $\mathcal{T}_{\rm hf}^0$ reduces the computational gain. In the numerical experiments, we provide extensive investigations of the influence of the coarse mesh on performance.
\end{itemize}

\section{Accelerated construction of  Galerkin ROMs for unsteady problems}
\label{sec:unsteady_problems}
We extend the acceleration strategy of section \ref{sec:steady_problems} to time-marching ROMs for unsteady PDEs.
 In view of the application considered in the numerical experiments, we here focus on Galerkin ROMs for discrete problems of the 
 form \eqref{eq:abstract_unsteadyPDE} with internal variables.
Several authors have developed MOR techniques for this class of equations --- see,  e.g.,
\cite{casenave2020nonintrusive,farhat2014dimensional,hernandez2017dimensional,iollo2022adaptive} and the references therein.
Here, we focus on the formulation presented in 
 \cite{agouzal2023projection} for quasi-static problems: our approach relies on Galerkin projection and on empirical quadrature for hyper-reduction;  as for the steady case, we envision that the acceleration strategy discussed in this work is also relevant for other MOR methods. 

\subsection{Vanilla POD-greedy for Galerkin time-marching ROMs}
\label{sec:PODgreedy}
 We denote by 
$u_{\mu}\in \mathbb{R}^{d}$ the displacement field,
by 
 $\sigma_{\mu}\in \mathbb{R}^{d\times d}$ the Cauchy tensor, by
 $\varepsilon_{\mu} = \nabla_{\rm s} u_{\mu}$ the strain tensor with 
 $\nabla_{\rm s} \bullet = \frac{1}{2} (\nabla \bullet + \nabla \bullet ^\top )$; we further introduce the vector of internal variables
 $\gamma_{\mu}\in \mathbb{R}^{d_{\rm int}}$. Then, we 
 introduce the quasi-static equilibrium equations (completed with suitable boundary and initial conditions)
\begin{equation}
\label{eq:boring_case_unsteady}
\left\{
\begin{array}{ll}
\displaystyle{-\nabla \cdot \sigma_{\mu} = f}& {\rm in \;\; \Omega} \\[3mm]
\displaystyle{\sigma_{\mu}  = \mathcal{F}_{\mu}^{\sigma}(\varepsilon_{\mu}, \gamma_{\mu})}
& {\rm in \;\; \Omega} \\[3mm]
\displaystyle{\dot{\gamma}_{\mu}  = \mathcal{F}_{\mu}^{\gamma}(\varepsilon_{\mu}, \gamma_{\mu})}
& {\rm in \;\; \Omega} \\
\end{array}
\right.
\end{equation} 
 where $\mathcal{F}_{\mu}^{\sigma}: \mathbb{R}^{d\times d} \times \mathbb{R}^{d_{\rm int}} \to \mathbb{R}^{d\times d}$ and $\mathcal{F}_{\mu}^{\gamma}: \mathbb{R}^{d\times d} \times \mathbb{R}^{d_{\rm int}} \to \mathbb{R}^{d_{\rm int}}$ are suitable parametric functions that encode the material constitutive law --- note that 
the Newton's law 
\eqref{eq:boring_case_unsteady}$_1$
does not include the  inertial term; 
 the temporal evolution is hence entirely driven by the constitutive law \eqref{eq:boring_case_unsteady}$_3$. Equation \eqref{eq:boring_case_unsteady} is discretized using the FE method in space and a one-step finite difference (FD) method in time. Given the 
parameter $\mu\in \mathcal{P}$, the  
 FE space $\mathcal{X}_{\rm hf}$ and the time grid $\{ t^{(k)} \}_{k=1}^K$, we seek the sequence
$\mathbb{u}_{\mu} = 
\{  {u}_{\mu}^{(k)}  \}_{k=0}^K$ such that
\begin{equation}
\label{eq:unsteadyPDE_FEM_a}
\mathfrak{R}_{\mu}^{(k)}(u_{\mu}^{(k)}, 
 v) = \sum_{\ell=1}^{N_{\rm e}} 
r_{\ell,\mu}^{(k), \rm e} \left(  u_{\mu}^{(k)},  v  \right) \; + \;
\sum_{j\in \mathcal{I}_{\rm bnd}}  r_{j, \mu}^{(k), \rm f} \left(
u_{\mu}^{(k)}  , v 
 \right),
\end{equation}
where the elemental residuals satisfy
\begin{equation}
\label{eq:unsteadyPDE_FEM_b}
\begin{array}{rl}
\displaystyle{
r_{\ell,\mu}^{(k), \rm e} \left(  w,  v  \right) 
=
\int_{\texttt{D}_{\ell}} \left(
\mathcal{F}_{\mu}^{\sigma}( \nabla_{\rm s} w, \gamma_{\mu}^{(k)}(w)    ) : \nabla_{\rm s} v
-
f \cdot v
\right) \, dx},
&
\ell=1,\ldots,N_{\rm e},
\\[4mm]
\displaystyle{
{\rm with} \; 
 \gamma_{\mu}^{(k)}(w) =  \mathcal{F}_{\mu,\Delta t}^{\gamma} }
 & 
 \displaystyle{
(\nabla_{\rm s} w, \varepsilon_{\mu}^{(k-1)}, \gamma_{\mu}^{(k-1)}) },
  \\ 
\end{array}
\end{equation}
$\mathcal{I}_{\rm bnd}\subset \{1,\ldots,N_{\rm f}\}$ are the indices of the boundary facets and
  the  facet  residuals
$\{ r_{j, \mu}^{(k), \rm f} \}_{j\in  \mathcal{I}_{\rm bnd}}$
incorporate boundary conditions.
Note that $\mathcal{F}_{\mu,\Delta t}^{\gamma}$ is the FD approximation of the constitutive law
\eqref{eq:boring_case_unsteady}$_3$.
 
 Given the reduced space 
 $\mathcal{Z}={\rm span} \{ \zeta_i \}_{i=1}^n \subset \mathcal{X}_{\rm hf}$,
 the time-marching hyper-reduced Galerkin ROM of \eqref{eq:unsteadyPDE_FEM_a}-\eqref{eq:unsteadyPDE_FEM_b} reads as:
 given $\mu\in \mathcal{P}$,  
 find
 $\widehat{\mathbb{u}}_{\mu} = 
\{  \widehat{u}_{\mu}^{(k)}  \}_{k=0}^K$ such that
\begin{equation}
\label{eq:unsteadyROM}
\mathfrak{R}_{\mu}^{(k), \rm eq}(\widehat{u}_{\mu}^{(k)}, 
 v) = \sum_{\ell=1}^{N_{\rm e}}  \rho_{\ell}^{\rm eq,e}
r_{\ell,\mu}^{(k), \rm e} \left(  \widehat{u}_{\mu}^{(k)},  v  \right) \; + \;
\sum_{j\in \mathcal{I}_{\rm bnd}}
\rho_{j}^{\rm eq,f}
  r_{j, \mu}^{(k), \rm f} \left(
\widehat{u}_{\mu}^{(k)}  , v 
 \right),
 \quad
 \forall \, v\in \mathcal{Z},
\end{equation}
for  $k=1,\ldots,K$,
 where $\boldsymbol{\rho}^{\rm eq,e} \in \mathbb{R}_+^{N_{\rm e}}$ and
 $\boldsymbol{\rho}^{\rm eq,f} \in \mathbb{R}_+^{N_{\rm f}}$ are suitably-chosen sparse vectors of weights.
 We observe that the Galerkin ROM \eqref{eq:unsteadyROM} does not reduce the total number of time steps: solution to \eqref{eq:unsteadyROM} hence requires the solution to a sequence of $K$ nonlinear problems of size $n$ and is thus likely much more expensive than the solution to (Petrov-) Galerkin ROMs for steady-state problems. Furthermore, several independent studies have shown that residual-based \emph{a posteriori}  error estimators for \eqref{eq:unsteadyROM} are typically much less sharp and reliable than their counterparts for steady-state problems. These observations have motivated the development of space-time formulations  \cite{urban2014improved}: to our knowledge, however, space-time methods have not been extended to problems with internal variables.

The abstract algorithm \ref{alg:weak_greedy} can be readily extended to time-marching ROMs for unsteady  PDEs: the POD-Greedy method \cite{haasdonk2008reduced} combines a greedy search in the parameter domain with a temporal compression based on POD. At each iteration $it$ of the algorithm, we update the reduced space $\mathcal{Z}$ using the newly-computed trajectory $\{  u_{\mu^{\star, it}}^{(k)}  \}_{k=1}^K$ and we update the quadrature rule. The reduced space $\mathcal{Z}$ can be updated using HAPOD (cf. section \ref{sec:steady_problems}) or using a hierarchical approach (cf. 
\cite[section 3.5]{haasdonk2017reduced}) that generates nested spaces: we refer to \cite[section 3.2.1]{iollo2022adaptive} for further details and extensive numerical investigations for a model problem with internal variables.
Here, we consider the nested approach of 
\cite{haasdonk2017reduced}: we denote by $n_{it}$ the dimension of the reduced space $\mathcal{Z}_{it}$ at the $it$-th iteration, and we denote by 
$n_{\rm new} = n_{it} - n_{it-1}$ the number of modes added at each iteration,
\begin{subequations} 
\label{eq:nested_space}
\begin{equation}
\mathcal{Z}_{it} = 
\mathcal{Z}_{it-1} \oplus   \mathcal{Z}^{\rm new},
 \quad
 {\rm where}
 \;
\mathcal{Z}^{\rm new} = 
\texttt{POD} 
\left( 
\{ \Pi_{\mathcal{Z}_{it-1}^\perp} u_{\mu^{\star, it}}^{(k)}  \}_{k=1}^K,  
n_{\rm new}(tol) ,
\; \;
(\cdot, \cdot)
\right),
\end{equation}
where $n_{\rm new} (tol)$ satisfies
$n_{\rm new}(tol) = $ 
\begin{equation}
\min \left\{
n' \, : \,
\frac{    
\sum_{k=1}^K \;
\Delta t^{(k)}
\| 
  u_{\mu^{\star, it}}^{(k)} - 
 \Pi_{  \mathcal{Z}_{it-1} \oplus   \mathcal{Z}_{n'}^{\rm new}  }    u_{\mu^{\star, it}}^{(k)} \|^2
}{\sum_{k=1}^K \;
\Delta t^{(k)}
\|   u_{\mu^{\star, it}}^{(k)} \|^2}
\leq
tol^2,
\; \;
\mathcal{Z}_{n'}^{\rm new} = 
{\rm span} \left\{
\zeta_{it,i}^{\rm new}
\right\}_{i=1}^{n'}
\right\},
\end{equation}
for some tolerance $tol>0$, with $\Delta t^{(k)}  =   t^{(k)}  - t^{(k-1)} $, $k=1,\ldots,K$. 
\end{subequations}
The quadrature rule is obtained using the same procedure described in section \ref{sec:progressive_EQETS_b}: exploiting \eqref{eq:nested_space}, it is easy to verify that the matrix $\mathbf{G}$  
 (cf. \eqref{eq:EQ_explained_c}) can be rewritten as (we omit the details)
\begin{equation}
\label{eq:matrixG_unsteady}
\mathbf{G}^{(it)}   = 
\left[
\begin{array}{lll}
\mathbf{G}^{\rm cnst,e}  &   & 0 \\ 
0 &   & \mathbf{G}^{\rm cnst,f}  \\ 
   &  \mathbf{G}_{\rm acc}^{(it-1)}  & \\[3mm]
   &  \mathbf{G}_{\rm new}^{(it)}  & \\[3mm]
\end{array}
\right],
{\rm with} \;\;
\left\{
\begin{array}{l}
\displaystyle{
\mathbf{G}_{\rm acc}^{(it-1)} \in \mathbb{R}^{(n_{it-1} (it - 1)K     ) \times (N_{\rm e}+N_{\rm f})     },
}\\[3mm]
\displaystyle{
\mathbf{G}_{\rm new}^{(it)}  \in \mathbb{R}^{(n_{\rm new} (it - 1)K + n_{it} K     ) \times (N_{\rm e}+N_{\rm f})     }.
}\\
\end{array}
\right.
\end{equation}
Note that $\mathbf{G}_{\rm acc}^{it}$ corresponds to the columns of the matrix $\mathbf{G}$ associated with the manifold accuracy constraints
(cf. \eqref{eq:accuracy_constraint}) at the $it$-th iteration of the greedy procedure.

\subsection{Acceleration strategy}
\label{sec:acceleration_unsteady} 
 
Acceleration of the POD-greedy algorithm for the Galerkin ROM \eqref{eq:unsteadyROM} relies on two independent building blocks: first,  the incremental construction of the quadrature rule; second,  a two-fidelity sampling strategy.
Re the quadrature procedure, we notice that
at each greedy iteration, 
 the matrix $\mathbf{G}$ in \eqref{eq:matrixG_unsteady}  
admits the decomposition in \eqref{eq:incrementalEQ_mat}. If the number of new modes is modest compared to $n_{it-1}$ the rows of 
$\mathbf{G}_{\rm new}^{(it)}$ are significantly less numerous than the columns of $\mathbf{G}^{(it)}$:  we can hence reduce the cost of the construction of $\mathbf{G}^{(it)}$ by keeping in memory the EQ matrix from the previous iteration; furthermore, the  decomposition
\eqref{eq:incrementalEQ_mat}
motivates the use of the active set of weights from the previous iterations to initialize the NNLS algorithm at the current iteration.
On the other hand, 
the extension of the  two-fidelity sampling strategy in Algorithm 
 \ref{alg:multifidelity_greedy} simply relies on the generalization of the strong -greedy procedure \ref{alg:strong_greedy} to unsteady problems, which is illustrated in the algorithm below.
 
\begin{algorithm}[H]                      
\caption{: POD-strong-greedy algorithm. }     
\label{alg:pod_strong_greedy}     

 \small
\begin{flushleft}
\emph{Inputs:}  $\{ 
\mathbb{u}_i =    \{  u_i^{(k)} \}_{k=1}^K 
  \}_{i=1}^{n_{\rm train}}$ snapshot set, 
$\texttt{maxit}$ maximum number of iterations,
$tol$ tolerance for data compression,
$(\cdot, \cdot)_{\rm sg}$ inner product.
\smallskip

\emph{Outputs:} 
$\mathcal{I}_{\rm sg}$ indices of selected snapshots.

\end{flushleft}                      

 \normalsize 

\begin{algorithmic}[1]
\State
Choose 
$\mathcal{Z}_0 = \emptyset$, 
$\mathcal{I}_{\rm sg}= \emptyset$,
set $\mathcal{I} _{\rm train}=\{1,\ldots,n_{\rm train}\}$.
\vspace{3pt}

\For {$it= 1, \ldots, \texttt{maxit}$}

\State
Compute $i^{\star} = {\rm arg} \max_{i\in \mathcal{I} _{\rm train}  } 
\sum_{k=1}^K\left( 
\min_{w\in \mathcal{Z}} 
\|  w - u_i^{(k)} \|_{\rm sg}^2 \right)$.
\vspace{3pt}

\State
Update the reduced space
$\mathcal{Z}_{it} = \texttt{data-compression}\left(
\mathcal{Z}_{it-1}, 
\{  u_{i^{\star}}^{(k)} \}_{k=1}^K, 
tol
\right)$.

\hfill  (cf. \eqref{eq:nested_space})
\vspace{3pt}

\State
$\mathcal{I}_{\rm sg}=  \mathcal{I}_{\rm sg} \cup   \{ i^{\star} \}$.
\vspace{3pt}
 
\EndFor
\end{algorithmic}
\bigskip

\end{algorithm}
  
Algorithm \ref{alg:pod_strong_greedy} takes as input a set of trajectories and returns the indices of the selected parameters. To shorten the notation, we here assume that the time grid is the same for all parameters: however, the algorithm can be readily extended to cope with parameter-dependent temporal discretizations. We further observe that the procedure depends on the data compression strategy employed to update the reduced space at each greedy iteration: in this work, we consider the same hierarchical strategy
 (cf. 
\cite[section 3.5]{haasdonk2017reduced}) employed in the POD-(weak-)greedy algorithm. 
  
\section{Numerical results}
\label{sec:numerics}

We present numerical results for 
a steady compressible inviscid flow past a LS89 blade (cf. section \ref{sec:LS89}), and for an unsteady nonlinear mechanics problem that simulates the 
long-time mechanical response of the standard section of a containment building under external loading (cf. section \ref{sec:NCB}).
Simulations of section \ref{sec:LS89}
are performed in Matlab 2022a \cite{MATLAB:2022} based on an in-house code, and executed over a commodity Linux workstation (RAM 32 GB, Intel i7 CPU 3.20 GHz x 12).
The HF simulations of section \ref{sec:NCB}
are performed using the FE software \texttt{code$\_$aster}
\cite{aster} and executed over a commodity 
Linux  workstation (RAM 32 GB, Intel i7-9850H CPU 2.60 GHz x 12);
on the other hand, 
the MOR procedure relies on an in-house \texttt{Python} code and is 
executed on a Windows  workstation
 (RAM 16 GB, Intel i7-9750H CPU 2.60 GHz x 12). 
 
\subsection{Transonic compressible flow past an LS89 blade}
\label{sec:LS89}

\subsubsection{Model problem}
\begin{figure}[h!]
\centering
\subfloat[]{ 
\includegraphics[width=.4\textwidth]{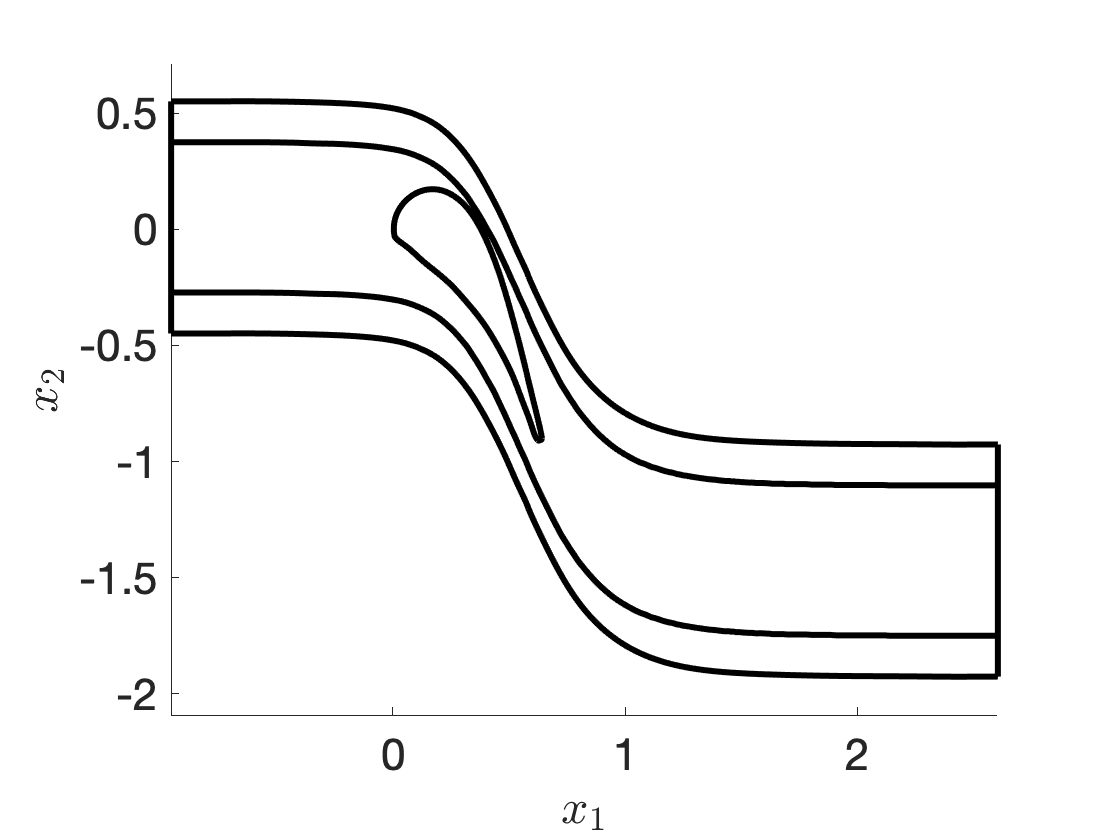}
}

\subfloat[]{ 
\includegraphics[width=.4\textwidth]{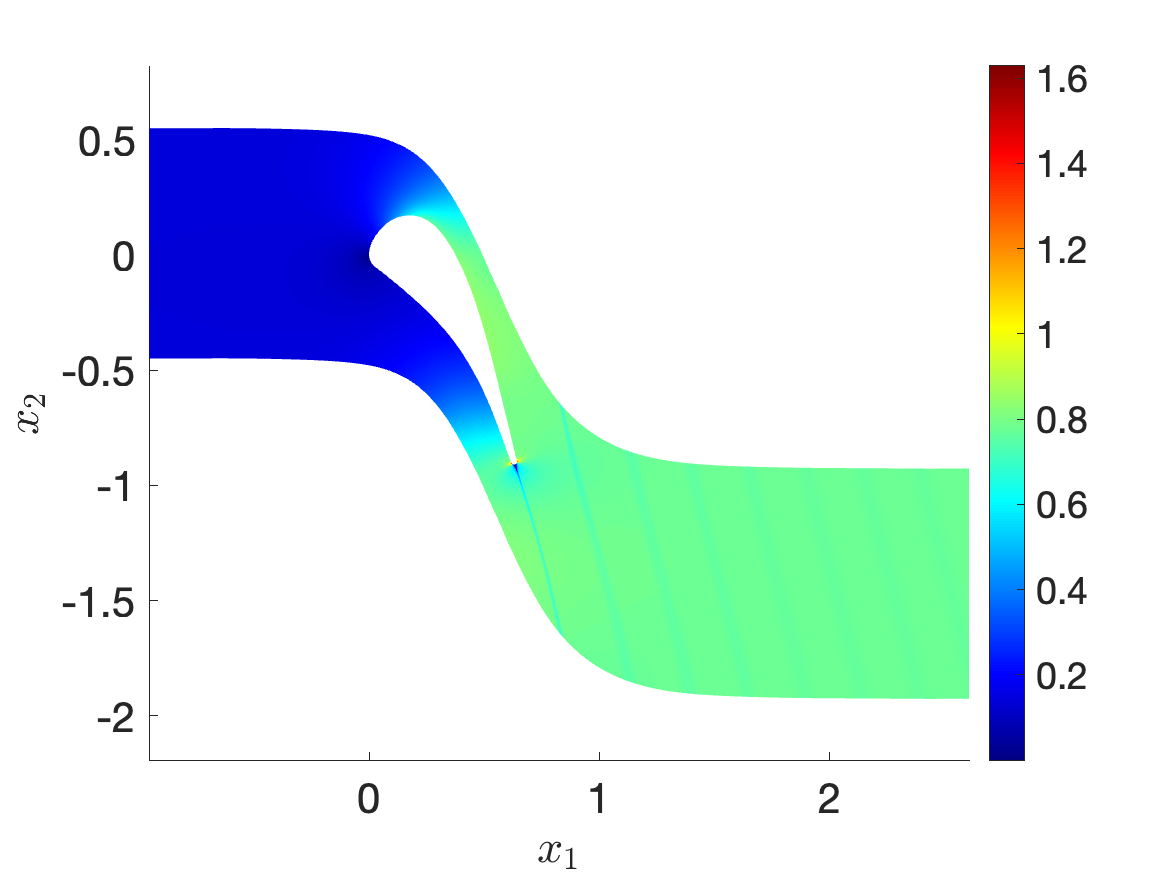}
}
~~
\subfloat[]{ 
\includegraphics[width=.4\textwidth]{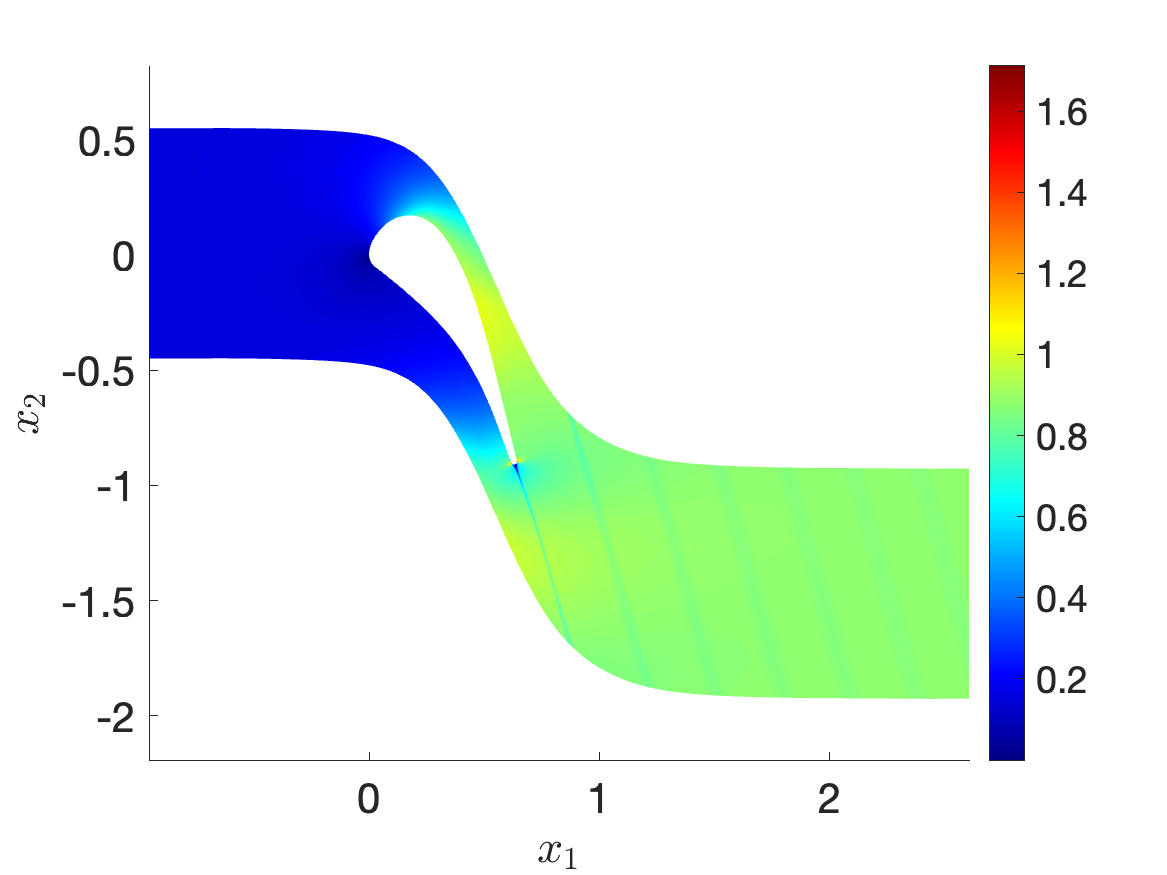}
}

\caption{inviscid flow past an array of LS89 turbine blades.
(a) partition associated with the geometric map.
(b)-(c) behavior of the Mach number for two parameter values.
}
\label{fig:visTB}
\end{figure}

We consider the problem of estimating the solution to the two-dimensional Euler equations past an array of LS89 turbine blades;
the same model problem is considered in \cite{taddei2023compositional}, for a different parameter range. 
We consider the computational domain depicted in Figure 
\ref{fig:visTB}(a); we prescribe total temperature, total pressure and flow direction at the inflow, static pressure at the outflow, non-penetration (wall) condition on the blade and periodic boundary conditions on the lower and upper  boundaries.
We study the sensitivity of the solution with respect to two parameters: the free-stream Mach number ${\rm Ma}_{\infty}$ and the height of the channel $H$, $\mu=[H,{\rm Ma}_{\infty}]$. We consider the parameter domain 
$\mathcal{P}=[0.9,1.1]\times [0.2,0.9]$. 

We deal with geometry variations through a piecewise-smooth mapping associated with the partition in 
Figure \ref{fig:visTB}(a).
We set $H_{\rm ref}=1$ and we define the curve $x_1\mapsto f_{\rm btm}(x_1)$ that describes the lower boundary $\Gamma_{\rm btm}$ of the domain  $\Omega=\Omega(H=1)$; then, we define $\widetilde{H}>0$ such that 
$x_1\mapsto f_{\rm btm}(x_1)+\widetilde{H}$ and
$x_1\mapsto f_{\rm btm}(x_1)+H-\widetilde{H}$ do not intersect the blade for any $H\in [0.9,1.1]$; finally, we define the geometric mapping
\begin{subequations}
\label{eq:geometric_mapping_turbine}
\begin{equation}
\Psi_H^{\rm geo}(x=[x_1,x_2])
=
\left[
\begin{array}{l}
x_1 \\
\psi_H^{\rm geo}(x)
\\
\end{array}
\right],
\end{equation}
where 
\begin{equation}
\psi_H^{\rm geo}(x)
=
\left\{
\begin{array}{ll}
o_1(x_1)+ C(H) \left(x_2 - o_1(x_1) \right) & x_2 < o_1(x_1),\\
o_2(x_1)+ C(H) \left(x_2 - o_2(x_1) \right) & x_2 > o_2(x_1),\\
x_2 & {\rm otherwise},
\\
\end{array}
\right.
\end{equation}
with 
$o_1(x_1)=f_{\rm btm}(x_1)+\widetilde{H}$,
$o_2(x_1)=f_{\rm btm}(x_1)+H_{\rm ref} - \widetilde{H}$ and
$C(H) = \frac{H-H_{\rm ref}}{2\widetilde{H}}+1$.
\end{subequations}
Figures \ref{fig:visTB}(b) and (c) show  the distribution of the Mach field for $\mu^{(1)}=[0.95,0.78]$ and 
$\mu^{(2)}=[1.05,0.88]$.
We notice that for large values of the free-stream Mach number the solution develops a normal shock on the upper side of the blade and two shocks at the trailing edge: effective approximation for higher values of the Mach number requires the use of nonlinear approximations and is beyond the scope of the present work.

\subsubsection{Results}

\begin{figure}[h!]
\centering
 \subfloat[mesh 1] 
{  \includegraphics[width=0.4\textwidth]
 {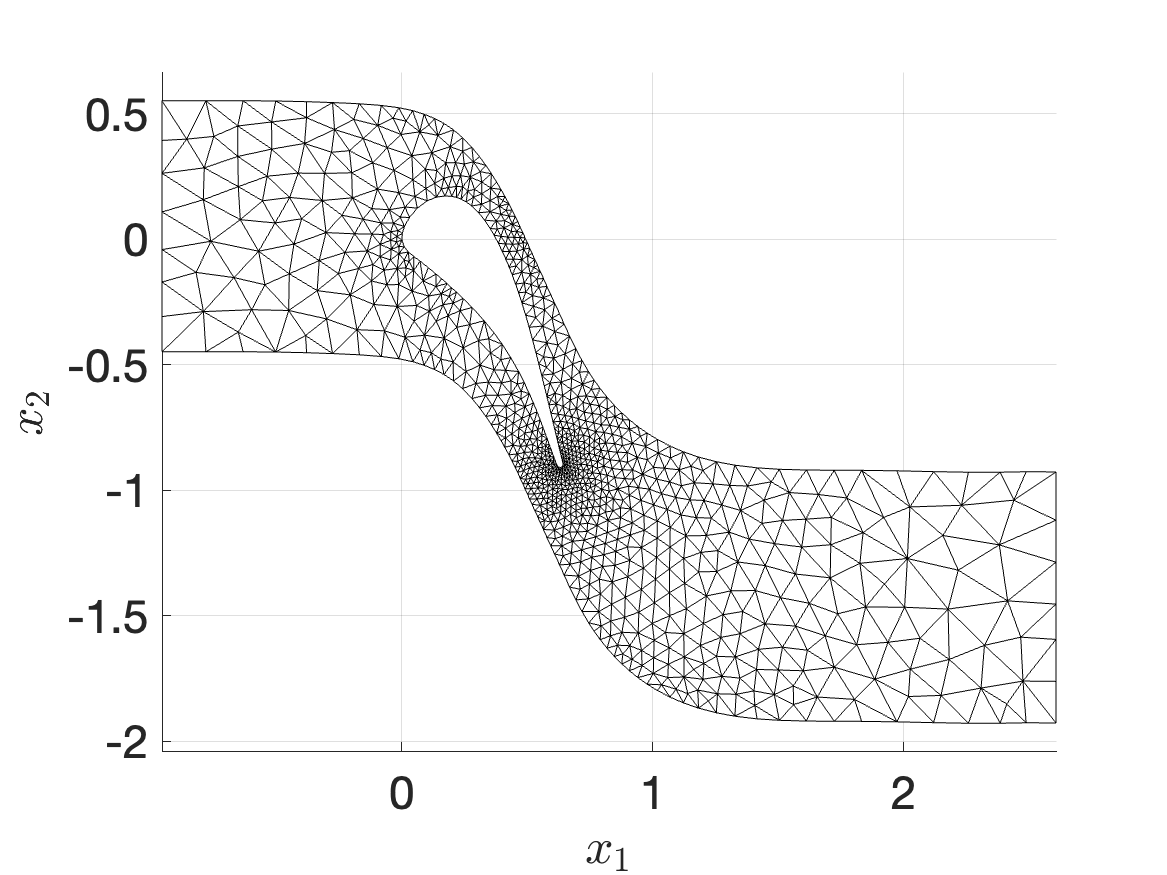}}
    ~~
 \subfloat[mesh 3] 
{  \includegraphics[width=0.4\textwidth]
 {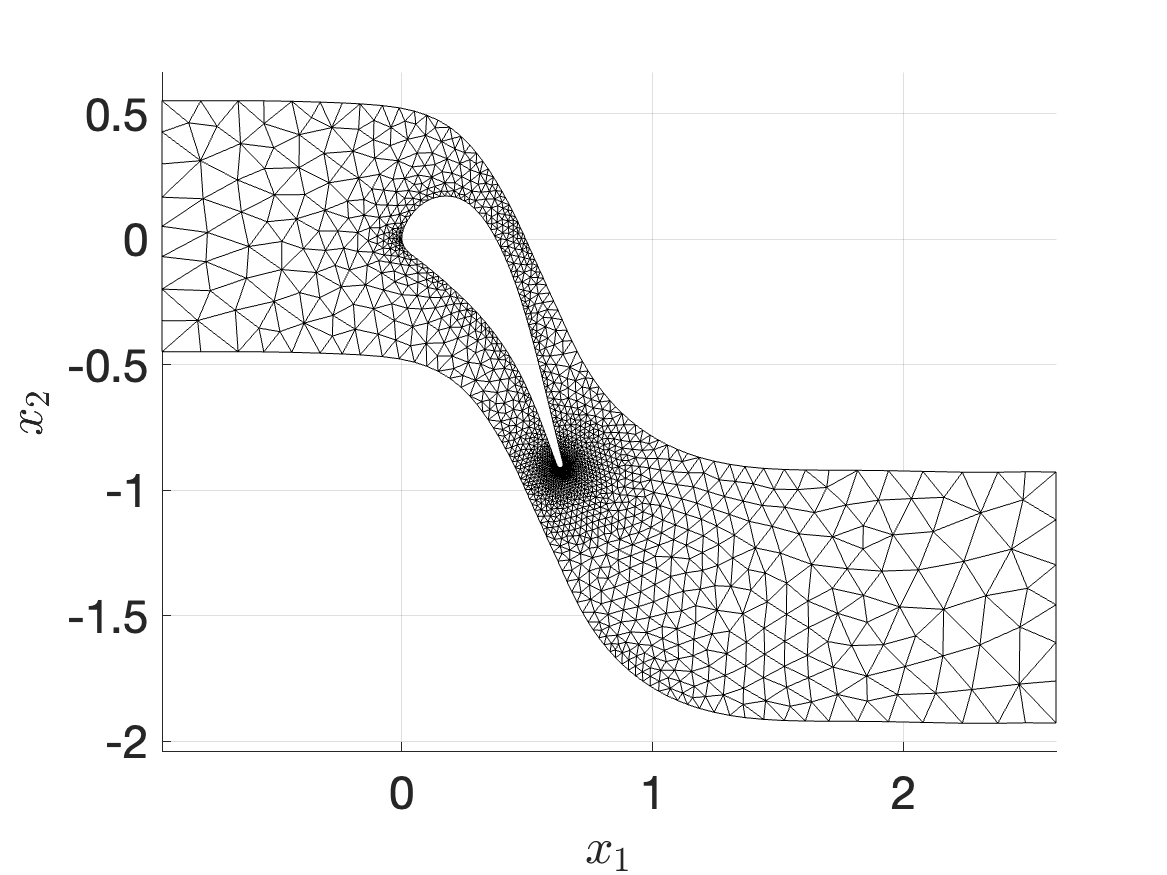}}

 \subfloat[mesh 6] 
{  \includegraphics[width=0.4\textwidth]
 {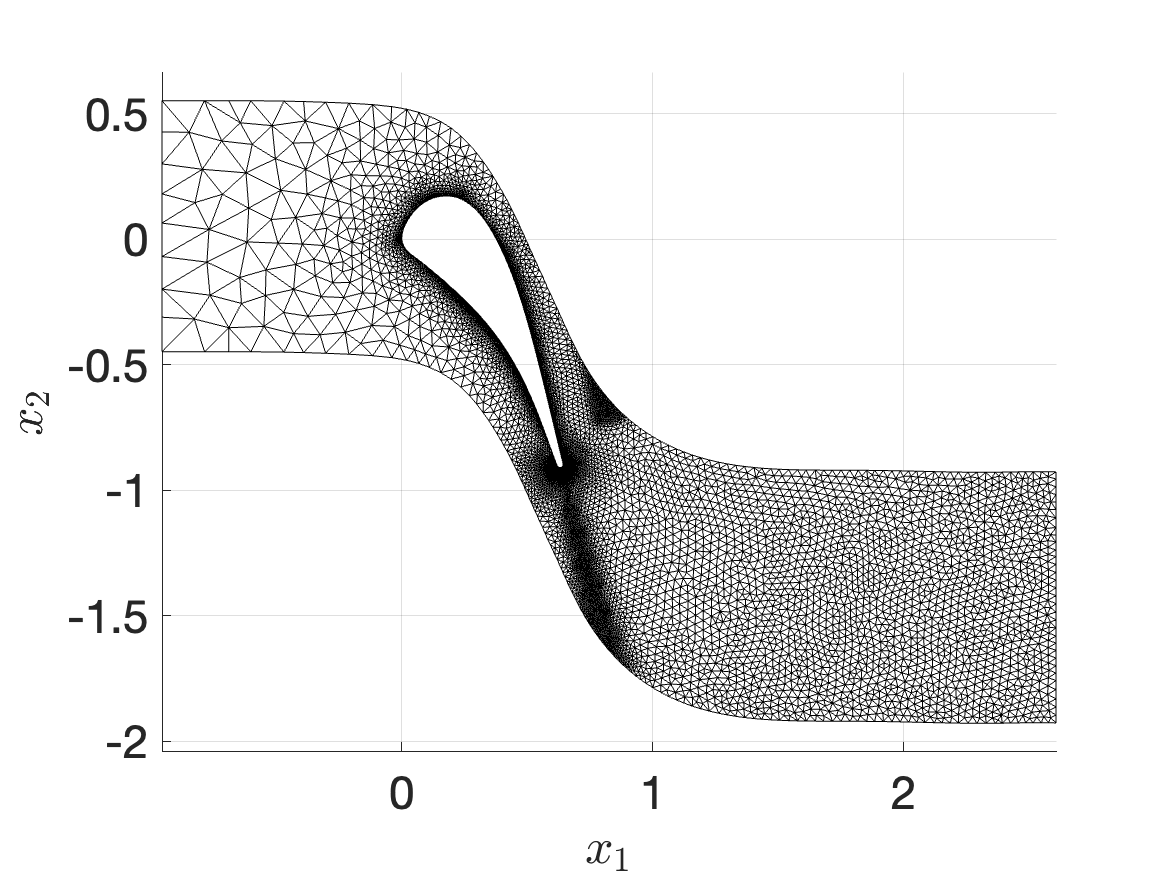}}
  ~~
\subfloat[]{
\begin{tikzpicture}[scale=0.5]
\begin{axis}[
xlabel = {\LARGE {$i$}},
  ylabel = {\LARGE {$E_{\mu}^{(i)}$}},
legend entries = {avg, max},
  line width=1.2pt,
  mark size=3.0pt,
ylabel style = {font=\large,yshift=1ex},
 yticklabel style = {font=\large,xshift=0ex},
xticklabel style = {font=\large,yshift=0ex},
legend style={at={(0.6,0.8)},anchor=west,font=\Large}
  ]
 
\addplot[color=red,mark=square]  table {data/euler/mesh/mean_error.dat};
%\addlegendentry{model}
  
\addplot[color=blue,mark=triangle*] table {data/euler/mesh/max_error.dat};
  
\end{axis}
\end{tikzpicture}
}
\caption{compressible flow past a LS89 blade.
(a)-(b)-(c) three computational meshes.
(d) behavior of the average and maximum error
\eqref{eq:grid_error}.
}
\label{fig:euler_vis_mesh}
\end{figure}
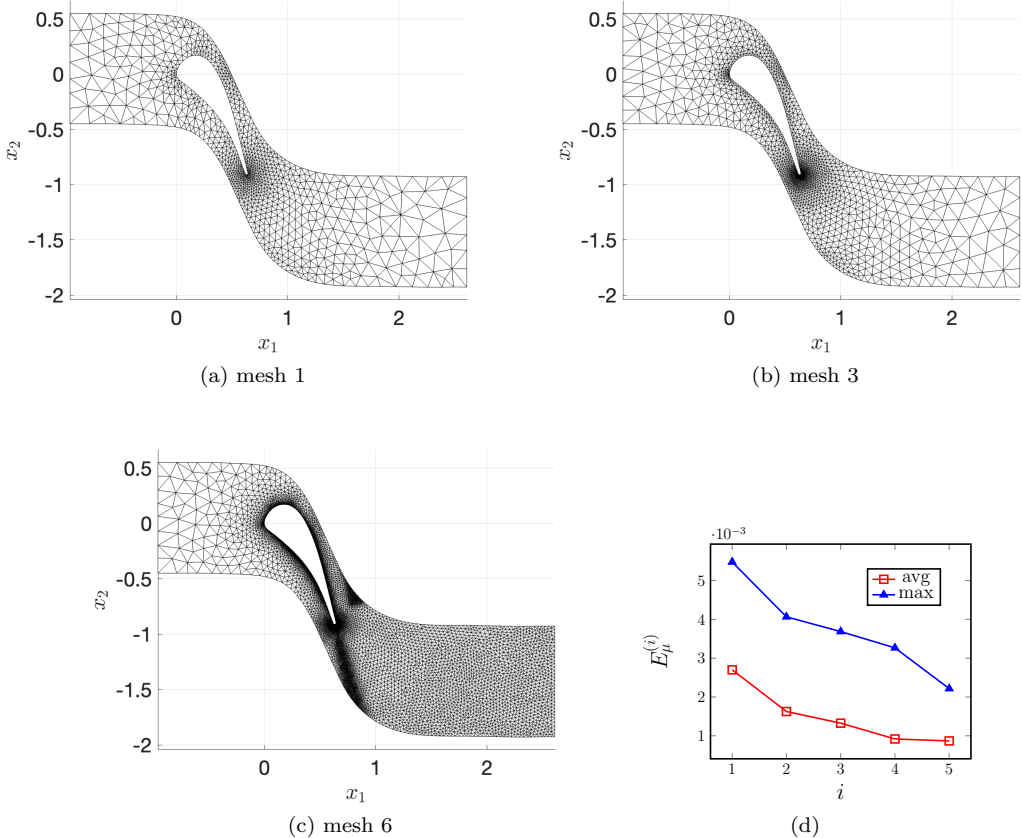

We consider a hierarchy of six P2 meshes with 
$N_{\rm e} = 1827, 2591, 3304, 4249,  7467, 16353$
elements, respectively;
Figures \ref{fig:euler_vis_mesh}(a)-(b)-(c) show three computational meshes. Figure \ref{fig:euler_vis_mesh}(d) show   maximum and mean errors between the HF solution $u_\mu^{{\rm hf}, (6)}$ and the corresponding HF solution associated with the $i$-th mesh,
\begin{equation}
\label{eq:grid_error}
E_{\mu}^{(i)} = \frac{\| u_\mu^{{\rm hf}, (i)} - u_\mu^{{\rm hf}, (6)}   \|}{\|  u_\mu^{{\rm hf}, (6)}   \|},
\end{equation}
over five randomly-chosen parameters in $\mathcal{P}$ . 
We remark that the error is measured in the reference configuration $\Omega$, which corresponds to $H=1$.
We consider the $L^2(\Omega)$ norm
$\|  \bullet \| = \sqrt{\int_{\Omega}  (\bullet)^2 \, dx }$.

%visualization of the meshes
   
\begin{figure}[h!]
\centering

%nbr its
\subfloat[]{
\begin{tikzpicture}[scale=0.5]
\begin{semilogyaxis}[
xlabel = {\LARGE {$n$}},
  ylabel = {\LARGE {nbr its EQP}},
legend entries = {incr, std},
  line width=1.2pt,
  mark size=3.0pt,
ylabel style = {font=\large,yshift=1ex},
 yticklabel style = {font=\large,xshift=0ex},
xticklabel style = {font=\large,yshift=0ex},
legend style={at={(0.01,0.6)},anchor=west,font=\Large}
  ]
 
\addplot[color=red,mark=square]  table {data/euler/it1/nbr_its.dat};
%\addlegendentry{model}
  
\addplot[color=blue,mark=triangle*] table {data/euler/it1/nbr_its_noincr.dat};
  
\end{semilogyaxis}
\end{tikzpicture}
}
~~~
%cost EQP
\subfloat[]{
\begin{tikzpicture}[scale=0.5]
\begin{axis}[
xlabel = {\LARGE {$n$}},
  ylabel = {\LARGE {cost EQP [s]}},
legend entries = {incr, std},
  line width=1.2pt,
  mark size=3.0pt,
ylabel style = {font=\large,yshift=1ex},
 yticklabel style = {font=\large,xshift=0ex},
xticklabel style = {font=\large,yshift=0ex},
legend style={at={(0.01,0.6)},anchor=west,font=\Large}
  ]
 
\addplot[color=red,mark=square]  table {data/euler/it1/costEQP.dat};
%\addlegendentry{model}
  
\addplot[color=blue,mark=triangle*] table {data/euler/it1/costEQP_noincr.dat};
  
\end{axis}
\end{tikzpicture}
}
~~~
%sampled elements
\subfloat[]{
\begin{tikzpicture}[scale=0.5]
\begin{axis}[
xlabel = {\LARGE {$n$}},
  ylabel = {\LARGE {$\%$ sampled weights}},
legend entries = {incr, std},
  line width=1.2pt,
  mark size=3.0pt,
ylabel style = {font=\large,yshift=1ex},
 yticklabel style = {font=\large,xshift=0ex},
xticklabel style = {font=\large,yshift=0ex},
legend style={at={(0.01,0.6)},anchor=west,font=\Large}
  ]
 
\addplot[color=red,mark=square]  table {data/euler/it1/nnz_tot.dat};
\addplot[color=blue,mark=triangle*] table {data/euler/it1/nnz_tot_noincr.dat};
\end{axis}
\end{tikzpicture}
}

%cost ES
\subfloat[]{
\begin{tikzpicture}[scale=0.5]
\begin{axis}[
xlabel = {\LARGE {$n$}},
  ylabel = {\LARGE {cost ES [s]}},
legend entries = {incr, std},
  line width=1.2pt,
  mark size=3.0pt,
ylabel style = {font=\large,yshift=1ex},
 yticklabel style = {font=\large,xshift=0ex},
xticklabel style = {font=\large,yshift=0ex},
legend style={at={(0.01,0.6)},anchor=west,font=\Large}
  ]
 
\addplot[color=red,mark=square]  table {data/euler/it1/costES.dat};

\addplot[color=blue,mark=triangle*] table {data/euler/it1/costES_noincr.dat};
 
\end{axis}
\end{tikzpicture}
}
~~~
 %L2err
\subfloat[]{
\begin{tikzpicture}[scale=0.5]
\begin{semilogyaxis}[
xlabel = {\LARGE {$n$}},
  ylabel = {\LARGE {rel $L^2$}},
legend entries = {incr, std},
  line width=1.2pt,
  mark size=3.0pt,
ymin=0.0001,   ymax=0.1,
ylabel style = {font=\large,yshift=1ex},
 yticklabel style = {font=\large,xshift=0ex},
xticklabel style = {font=\large,yshift=0ex},
legend style={at={(0.6,0.6)},anchor=west,font=\Large}
  ]
 
\addplot[color=red,mark=square]  table {data/euler/it1/L2err.dat};

\addplot[color=blue,mark=triangle*] table {data/euler/it1/L2err_noincr.dat};
 
\end{semilogyaxis}
\end{tikzpicture}
}
~~~
 %comp cost
\subfloat[]{
\begin{tikzpicture}[scale=0.5]
\begin{axis}[
xlabel = {\LARGE {$n$}},
  ylabel = {\LARGE {online cost [s]}},
legend entries = {incr, std},
  line width=1.2pt,
  mark size=3.0pt,
ylabel style = {font=\large,yshift=1ex},
 yticklabel style = {font=\large,xshift=0ex},
xticklabel style = {font=\large,yshift=0ex},
legend style={at={(0.01,0.6)},anchor=west,font=\Large}
  ]
 
\addplot[color=red,mark=square]  table {data/euler/it1/comp_cost.dat};

\addplot[color=blue,mark=triangle*] table {data/euler/it1/comp_cost_noincr.dat};
 
\end{axis}
\end{tikzpicture}
}

\caption{compressible flow past a LS89 blade. Progressive construction of  quadrature  rule and test space, coarse mesh ($N_{\rm e} = 1827$).}
\label{fig:euler_EQP_mesh1}
\end{figure}
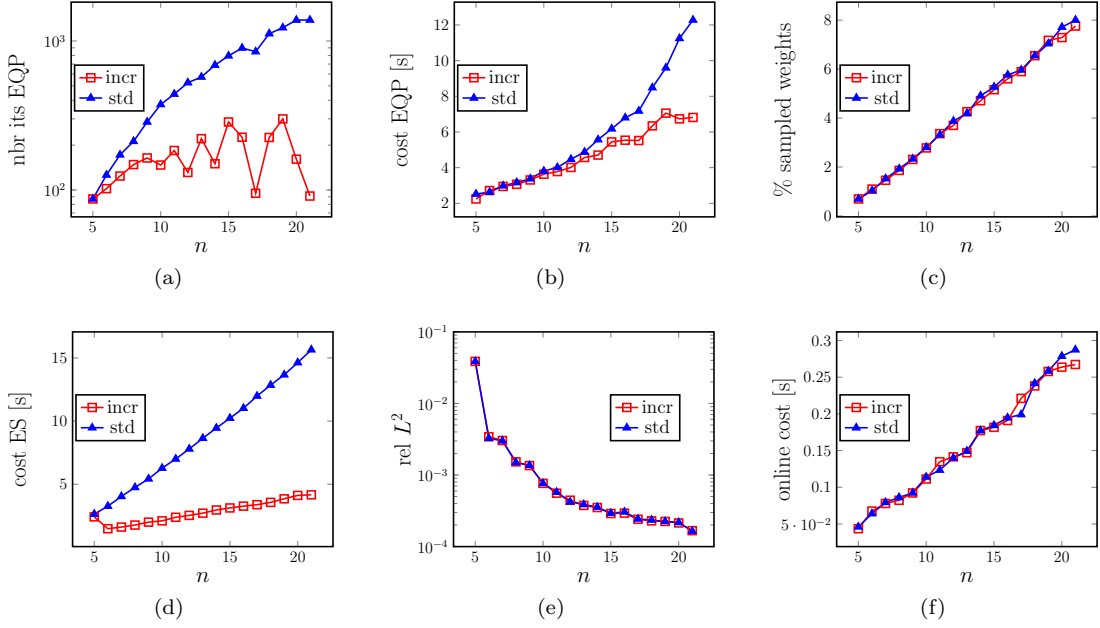

   Figure \ref{fig:euler_EQP_mesh1} compares the performance of the standard (``std'') greedy method with  the performance of the incremental (``incr'') procedures of section \ref{sec:progressive_EQETS} for the coarse mesh (mesh 1).
 In all the tests below, we consider a training space $\mathcal{P}_{\rm train}$ based on a ten by ten equispaced discretization of the parameter domain. 
 Figure \ref{fig:euler_EQP_mesh1}(a)  shows the number of iterations of the NNLS algorithm \ref{alg:nnls}, while 
   Figure \ref{fig:euler_EQP_mesh1}(b) shows the wall-clock cost (in seconds) on a commodity laptop, and
   Figure \ref{fig:euler_EQP_mesh1}(c)  shows  the percentage of sampled weights
   $\frac{\texttt{nnz}(\boldsymbol{\rho}^{\rm eq})}{N_{\rm e}+N_{\rm f}} \times 100 \%$.
   We observe that the proposed initialization of the active set method leads to a significant reduction of the total number of iterations, and to a non-negligible reduction of the total wall-clock cost\footnote{The computational cost reduction is not as significant  as the reduction in the total number of iterations, because the cost per iteration depends on the size of the least-square problem to be solved (cf. Line 11, Algorithm \ref{alg:nnls}), which increases as we increase the cardinality of the active set.}, without affecting the performance of the method. 
Figure \ref{fig:euler_EQP_mesh1}(d)  shows the cost of constructing the test space: we notice that the progressive construction of the test space enables a significant reduction of offline costs.  
 Finally, Figure \ref{fig:euler_EQP_mesh1}(e)   shows the averaged out-of-sample performance of the LSPG ROM over $n_{\rm test}=20$ randomly-chosen  parameters
$\mathcal{P}_{\rm test}$, 
$$
E_{\rm avg} = \frac{1}{n_{\rm test}} \sum_{\mu\in \mathcal{P}_{\rm test}}
\frac{\| u_\mu^{{\rm hf}, (1)} - \widehat{u}_\mu   \|}{\|  u_\mu^{{\rm hf}, (1)}   \|},
$$
and Figure \ref{fig:euler_EQP_mesh1}(f) shows  the online costs: we observe that the standard and the incremental approaches lead to nearly-equivalent results for all values of the ROB size $n$ considered.
   
   %performance  (mesh 1)  

\begin{figure}[h!]
\centering

%nbr its
\subfloat[]{
\begin{tikzpicture}[scale=0.5]
\begin{semilogyaxis}[
xlabel = {\LARGE {$n$}},
  ylabel = {\LARGE {nbr its EQP}},
legend entries = {incr, std},
  line width=1.2pt,
  mark size=3.0pt,
xmin=5,   xmax=22,
ylabel style = {font=\large,yshift=1ex},
 yticklabel style = {font=\large,xshift=0ex},
xticklabel style = {font=\large,yshift=0ex},
legend style={at={(0.01,0.6)},anchor=west,font=\Large}
  ]
 
\addplot[color=red,mark=square]  table {data/euler/it6/nbr_its.dat};
%\addlegendentry{model}
  
\addplot[color=blue,mark=triangle*] table {data/euler/it6/nbr_its_noincr.dat};
  
\end{semilogyaxis}
\end{tikzpicture}
}
~~~
%cost EQP
\subfloat[]{
\begin{tikzpicture}[scale=0.5]
\begin{axis}[
xlabel = {\LARGE {$n$}},
  ylabel = {\LARGE {cost EQP [s]}},
legend entries = {incr, std},
  line width=1.2pt,
  mark size=3.0pt,
xmin=5,   xmax=22,
ylabel style = {font=\large,yshift=1ex},
 yticklabel style = {font=\large,xshift=0ex},
xticklabel style = {font=\large,yshift=0ex},
legend style={at={(0.01,0.6)},anchor=west,font=\Large}
  ]
 
\addplot[color=red,mark=square]  table {data/euler/it6/costEQP.dat};
%\addlegendentry{model}
  
\addplot[color=blue,mark=triangle*] table {data/euler/it6/costEQP_noincr.dat};
  
\end{axis}
\end{tikzpicture}
}
~~~
%sampled elements
\subfloat[]{
\begin{tikzpicture}[scale=0.5]
\begin{axis}[
xlabel = {\LARGE {$n$}},
  ylabel = {\LARGE {$\%$ sampled weights}},
legend entries = {incr, std},
  line width=1.2pt,
  mark size=3.0pt,
xmin=5,   xmax=22,
ylabel style = {font=\large,yshift=1ex},
 yticklabel style = {font=\large,xshift=0ex},
xticklabel style = {font=\large,yshift=0ex},
legend style={at={(0.01,0.6)},anchor=west,font=\Large}
  ]
 
\addplot[color=red,mark=square]  table {data/euler/it6/nnz_tot.dat};
\addplot[color=blue,mark=triangle*] table {data/euler/it6/nnz_tot_noincr.dat};
\end{axis}
\end{tikzpicture}
}

%cost ES
\subfloat[]{
\begin{tikzpicture}[scale=0.5]
\begin{axis}[
xlabel = {\LARGE {$n$}},
  ylabel = {\LARGE {cost ES [s]}},
legend entries = {incr, std},
  line width=1.2pt,
  mark size=3.0pt,
xmin=5,   xmax=22,
ylabel style = {font=\large,yshift=1ex},
 yticklabel style = {font=\large,xshift=0ex},
xticklabel style = {font=\large,yshift=0ex},
legend style={at={(0.01,0.6)},anchor=west,font=\Large}
  ]
 
\addplot[color=red,mark=square]  table {data/euler/it6/costES.dat};

\addplot[color=blue,mark=triangle*] table {data/euler/it6/costES_noincr.dat};
 
\end{axis}
\end{tikzpicture}
}
~~~
 %L2err
\subfloat[]{
\begin{tikzpicture}[scale=0.5]
\begin{semilogyaxis}[
xlabel = {\LARGE {$n$}},
  ylabel = {\LARGE {rel $L^2$}},
legend entries = {incr, std},
  line width=1.2pt,
  mark size=3.0pt,
xmin=5,   xmax=22,
ymin=0.0001,   ymax=0.1,
ylabel style = {font=\large,yshift=1ex},
 yticklabel style = {font=\large,xshift=0ex},
xticklabel style = {font=\large,yshift=0ex},
legend style={at={(0.6,0.6)},anchor=west,font=\Large}
  ]
 
\addplot[color=red,mark=square]  table {data/euler/it6/L2err.dat};

\addplot[color=blue,mark=triangle*] table {data/euler/it6/L2err_noincr.dat};
 
\end{semilogyaxis}
\end{tikzpicture}
}
~~~
 %comp cost
\subfloat[]{
\begin{tikzpicture}[scale=0.5]
\begin{axis}[
xlabel = {\LARGE {$n$}},
  ylabel = {\LARGE {online cost [s]}},
legend entries = {incr, std},
  line width=1.2pt,
  mark size=3.0pt,
xmin=5,   xmax=22,
ylabel style = {font=\large,yshift=1ex},
 yticklabel style = {font=\large,xshift=0ex},
xticklabel style = {font=\large,yshift=0ex},
legend style={at={(0.01,0.6)},anchor=west,font=\Large}
  ]
 
\addplot[color=red,mark=square]  table {data/euler/it6/comp_cost.dat};

\addplot[color=blue,mark=triangle*] table {data/euler/it6/comp_cost_noincr.dat};
 
\end{axis}
\end{tikzpicture}
}

\caption{compressible flow past a LS89 blade. Progressive construction of  quadrature  rule and test space; fine mesh  ($N_{\rm e} = 16353$).}
\label{fig:euler_EQP_mesh6}
\end{figure}
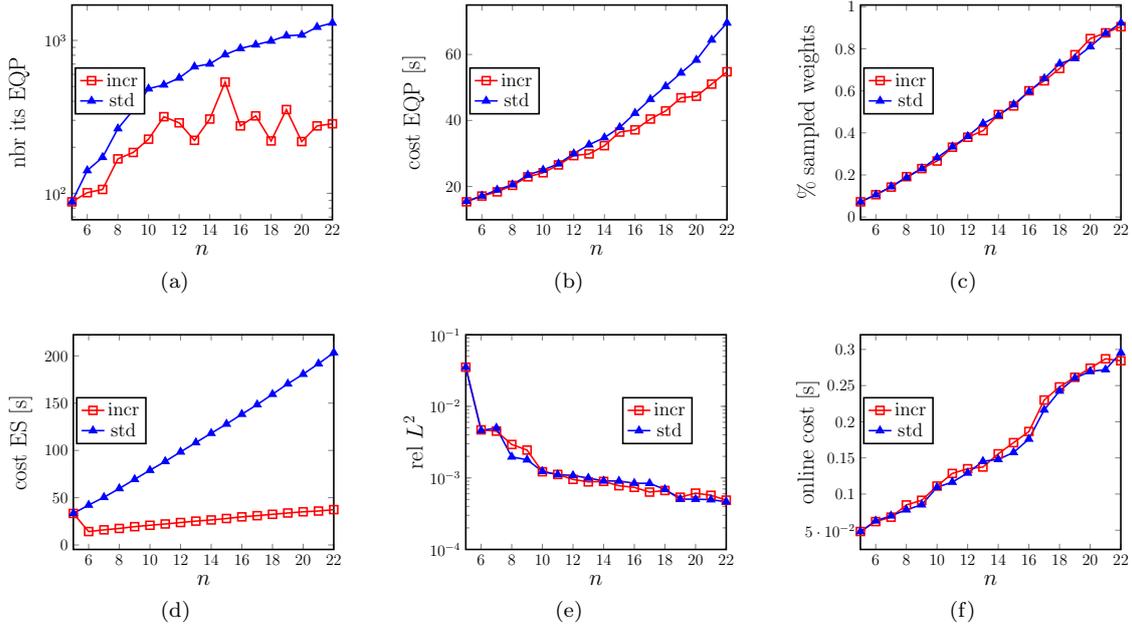

Figure \ref{fig:euler_EQP_mesh6}
replicates the tests of Figure \ref{fig:euler_EQP_mesh1} for the fine mesh  (mesh 6).
As for the coarser grid, we find that the progressive construction of the test space and of the quadrature rule do not hinder online performance of the LSPG ROM and ensure significant offline savings.
We notice that the relative error over the test set is significantly larger: reduction of the mesh size leads to more accurate approximations of  sharp   features
--- shocks, wakes --- that are troublesome for linear approximations.
On the other hand, we notice that online costs are nearly the same as for the coarse mesh for the corresponding value of the ROB size $n$: this proves the effectiveness of the hyper-reduction procedure.
 
%performance (mesh6)
   
\begin{figure}[h!]
\centering
 \subfloat[] 
{  \includegraphics[width=0.4\textwidth]
 {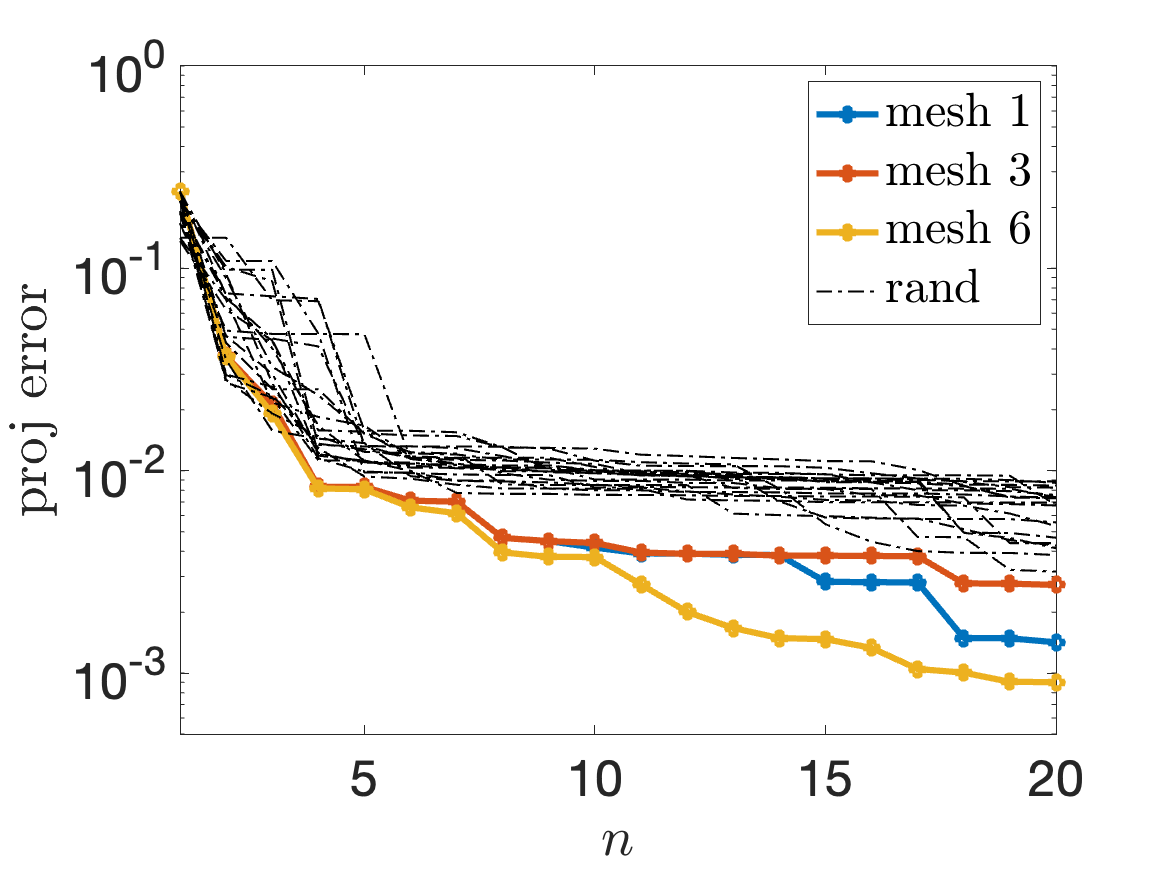}}

 \subfloat[mesh 1] 
{  \includegraphics[width=0.4\textwidth]
 {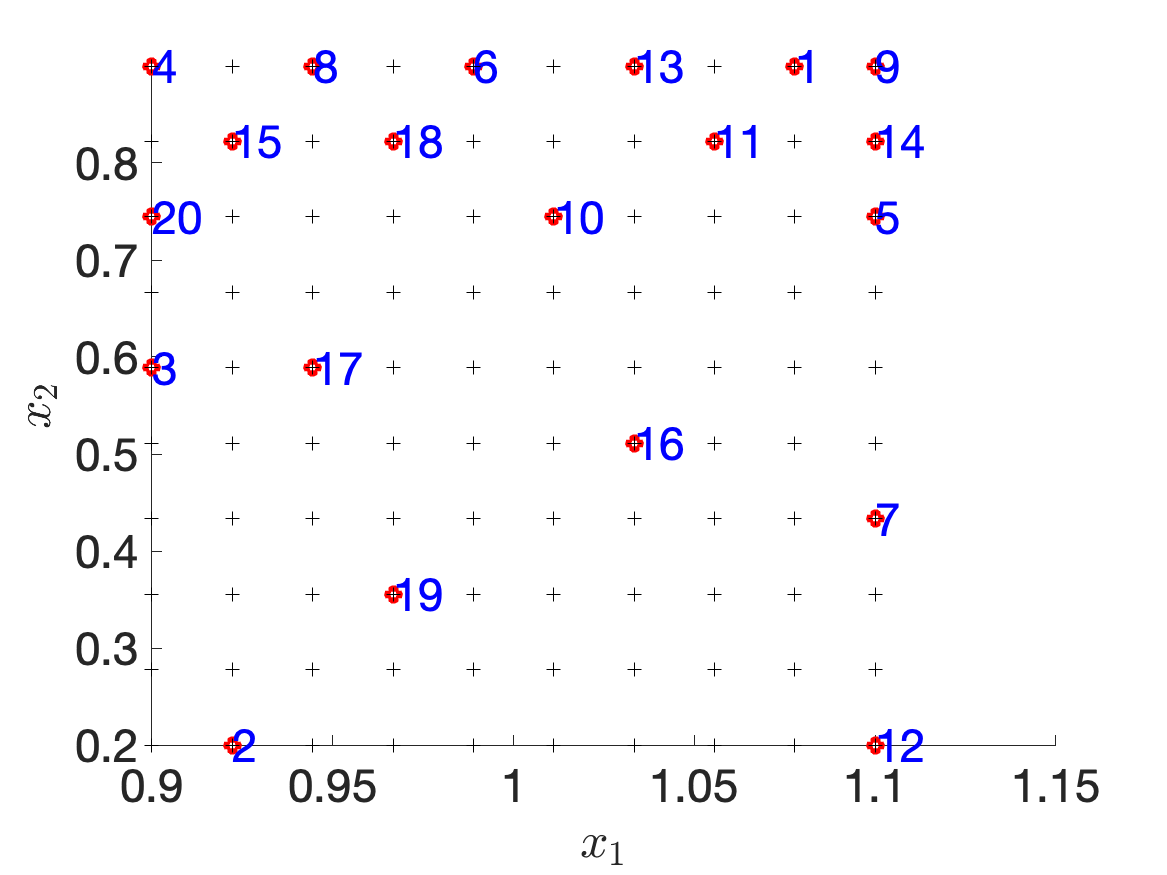}}
    ~~
 \subfloat[mesh 6] 
{  \includegraphics[width=0.4\textwidth]
 {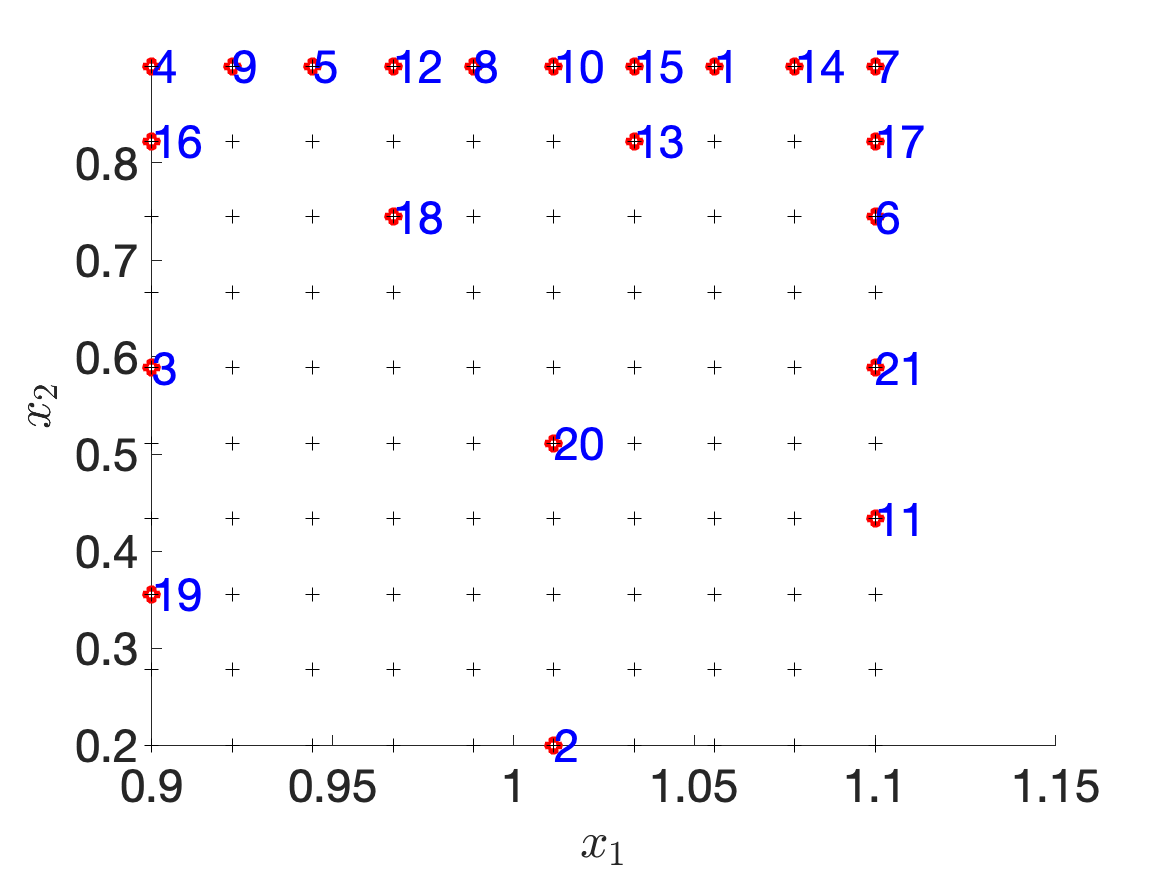}}

\caption{compressible flow past a LS89 blade; sampling.
(a)  behavior of $E_{n}^{{\rm proj},(i)}$ 
\eqref{eq:proj_error_sampling}
for six different choices of the coarse mesh, and for random samples.
(b)-(c)  parameters $\{ \mu^{\star,j} \}_j$ selected  by Algorithm \ref{alg:multifidelity_greedy} for two different coarse meshes.
}
\label{fig:euler_vis_sampling}
\end{figure}

Figure \ref{fig:euler_vis_sampling} investigates the effectiveness of the sampling strategy based on the strong-greedy algorithm.  We rely on Algorithm \ref{alg:multifidelity_greedy} to identify the training set of parameters $\mathcal{P}_\star$ for different choices of the coarse mesh $\mathcal{T}_{\rm hf}^0=\mathcal{T}_{\rm hf}^{(i)}$, for $i=1,\ldots,6$ and for 
 $\mathcal{T}_{\rm hf} =\mathcal{T}_{\rm hf}^{(6)}$.
 Then, we measure performance in terms of the maximum relative projection error over $\mathcal{P}_{\rm train}$,
\begin{equation}
\label{eq:proj_error_sampling}
E_{n}^{{\rm proj},(i)} = 
\max_{\mu\in \mathcal{P}_{\rm train}}
\inf_{\zeta\in \mathcal{Z}_n^{(i)}}
\frac{\| u_\mu^{{\rm hf}, (6)} - \zeta   \|}{\|  u_\mu^{{\rm hf}, (6)}  \|},
\quad
{\rm with}
\;\;
\mathcal{Z}_{n'}^{(i)}
=
{\rm span} \left\{
u_\mu^{{\rm hf}, (6)} : 
\mu\in \mathcal{P}_{\star,n'}^{(i)}
\right\},
\end{equation} 
where $\mathcal{P}_{\star,n'}^{(i)}$ is the set of the first $n'$ parameters selected through Algorithm \ref{alg:multifidelity_greedy} based on the coarse mesh  
$\mathcal{T}_{\rm hf}^{(i)}$.
Figure \ref{fig:euler_vis_sampling}(a) shows the behavior of the projection error
$E_{n}^{{\rm proj},(i)}$ \eqref{eq:proj_error_sampling} for three different choices of the coarse mesh; to provide a concrete reference, we also report the performance of twenty
sequences of reduced spaces obtained by randomly selecting sequences of parameters in $\mathcal{P}_{\rm train}$. 
 Figures \ref{fig:euler_vis_sampling}(b) and (c) show the parameters selected through Algorithm \ref{alg:multifidelity_greedy}  for two different choices of the coarse mesh: we observe that the selected parameters are clustered in the proximity of ${\rm Ma}_{\infty}=0.9$.
   
%sampling
   
Table \ref{tab:LS89} compares the costs of the standard weak-greedy algorithm (``vanilla''), 
 the  weak-greedy algorithm with progressive construction of the test space and the quadrature rule (``incr''), and the two-fidelity algorithm \ref{alg:multifidelity_greedy} with coarse mesh given by $\mathcal{T}_{\rm hf}^0= \mathcal{T}_{\rm hf}^{(1)}$ (``incr+MF'').
 To ensure a fair comparison, we impose that the final ROM has the same number of modes (twenty) for all cases.
Training of the coarse ROM
in Algorithm \ref{alg:multifidelity_greedy}  is based on the weak-greedy algorithm with progressive construction of the test space and the quadrature rule, with tolerance $\texttt{tol}=10^{-3}$: this leads to a coarse ROM with $n_0=14$ modes, which corresponds to an initial training set $\mathcal{P}_\star$ of cardinality $14$ in the greedy method (cf. Line 5, Algorithm \ref{alg:multifidelity_greedy}).
The ROMs associated with three different training strategies show comparable 
performance  in terms of online cost and $L^2$ errors.

\begin{table}[h!]
\centering
\begin{tabular}{|l|c|c|c|c|c|c|c|}
\hline 
\textbf{fine mesh $\mathcal{T}_{\rm hf}^{(5)}$}  
&  ROB & ES & EQP &  greedy search & HF solves & overhead & total \\[0.5ex]  \hline \hline 
vanilla &  $1.1$ &  $714.5$ &   $247.4$ &  $930.6$ &  $3554.8$ &  $0.6$ &  $5449.0$ \\[0.5ex]  \hline   
incr &  $1.1$ &  $177.9$ &   $221.2$ &  $928.4$ &  $3504.8$ &  $0.6$ &  $4834.0$ \\[0.5ex]  \hline   
incr+MF &  $0.9$ &  $131.0$ &   $121.9$ &  $362.7$ &  $2030.0$ &  $692.1$ &  $3338.5$ \\[0.5ex]  \hline   
\end{tabular}

\medskip

\begin{tabular}{|l|c|c|c|c|c|c|c|}
\hline 
\textbf{fine mesh $\mathcal{T}_{\rm hf}^{(6)}$}
& ROB & ES & EQP &  greedy search & HF solves & overhead & total \\[0.5ex]  \hline \hline 
vanilla &  $3.4$ &  $2067.8$ &   $668.4$ &  $2091.2$ &  $10997.4$ &  $2.3$ &  $15830.5$ \\[0.5ex]  \hline   
incr &  $3.2$ &  $444.8$ &   $538.4$ &  $1952.8$ &  $10529.1$ &  $2.3$ &  $13470.6$ \\[0.5ex]  \hline   
incr+MF &  $3.6$ &  $366.4$ &   $362.5$ &  $905.8$ &  $9597.7$ &  $699.9$ &  $11935.9$ \\[0.5ex]  
\hline   
\end{tabular}

\caption{compressible flow past a LS89 blade. Overview of offline costs for three computational strategies and two different fine meshes (coarse mesh: $\mathcal{T}_{\rm hf}^{(1)}$).}
\label{tab:LS89}
\end{table}

For the fine mesh $\mathcal{T}_{\rm hf}^{(6)}$, 
the two-fidelity training leads to a reduction of offline costs of roughly $25\%$ with respect to the vanilla implementation and of roughly  $10\%$ with respect to the incremental implementation;
for the fine mesh $\mathcal{T}_{\rm hf}^{(5)}$ (which has one half as many elements as the finer grid), 
the two-fidelity training leads to a reduction of offline costs of roughly $39\%$ with respect to the vanilla implementation and of roughly  
$31\%$ with respect to the incremental implementation.
In particular, we notice that the initialization based on the coarse model --- which is the same for both cases --- is significantly more effective for the HF model associated with mesh 5 than for the HF model associated with  mesh 6. The empirical finding strengthens the observation made in section \ref{sec:multifidelity_sampling} that the choice of the coarse approximation is a compromise between overhead costs --- which increase  as we increase  the size of the coarse mesh --- and accuracy of the coarse solution.

We remark that for the first two cases the HF solver is initialized using the solution for the closest parameter in the training set
$$
u_{\mu^{\star,n}}^0 = u_{\bar{\mu}}^{\rm hf},
\quad
{\rm with} \;\;
\bar{\mu}
=
{\rm arg}\min_{\mu'\in 
\{
\mu^{\star,i}
\}_{i=1}^{n-1}
 }
 \| \mu^{\star,n} - \mu' \|_2,
$$
for $n=2,3,\ldots$; 
for $n=1$, we rely on a coarse solver and on a continuation strategy with respect to the Mach number:
  this  initialization is inherently sequential.
On the other hand, in the two-fidelity procedure, the HF solver is initialized using the reduced-order solver that is trained using 
 coarse HF data: this choice enables 
 trivial parallelization of the HF solves for the initial parameter sample 
 $\mathcal{P}_\star$. 
We hence expect to achieve higher computational savings for parallel computations.

\subsection{Long-time mechanical response of the standard section of a containment building under external loading}
\label{sec:NCB}

\subsubsection{Model problem}
We study the long-time mechanical response of a three-dimensional standard section of a nuclear power plant containment building: the highly-nonlinear mechanical response is activated by thermal effects; its simulation requires the coupling of thermal, hydraulic   and mechanical (THM) responses.
A thorough presentation of the mathematical model and of the MOR procedure is provided in \cite{agouzal2024projection}.
The ultimate goal of the simulation  is to predict the temporal behavior of several quantities of interest (QoIs), such as water saturation in concrete, delayed deformations, and stresses: these QoIs are directly related to the leakage rate, whose estimate is of paramount importance for the design of NCBs.
The deformation field is also important to conduct validation and calibration studies 
against real-world data.

Following \cite{bouhjiti2018analyse}, we consider a weak THM coupling procedure
to model deferred deformations within the material; weak coupling is appropriate for large structures under normal operational loads. 
The MOR process is exclusively applied to estimate the mechanical response: the results from thermal and hydraulic calculations are indeed used as input data for  the mechanical calculations,
which constitute the computational bottleneck of the entire simulation workflow.
To model the mechanical response of the concrete structure, we  consider  a three-dimensional nonlinear rheological creep model with internal variables;
on the other hand, we consider a one-dimensional linear elastic model for the prestressing cables:
the state variables are hence the  displacement field of the three-dimensional concrete structure and of the one-dimensional steel cables.
We assume that the whole structure satisfies  the small-strain and small-displacement hypotheses.
To establish connectivity between concrete and steel nodes, a kinematic linkage is implemented: a point within the steel structure and its corresponding point within the concrete structure are assumed to share identical displacements.

\begin{figure}[h!]
\centering
%normal effort horizontal cable
\subfloat[]{
\begin{tikzpicture}[scale=0.5]
\begin{axis}[
	xlabel = {\LARGE {$t$ (s)}},
	ylabel = {\LARGE {$N_{\rm H2}$}},
	legend entries = {$\mu^{(1)}$, $\mu^{(2)}$, $\mu^{(3)}$},
	line width=1.2pt,
	mark size=3.0pt,
	ylabel style = {font=\large,yshift=1ex},
	yticklabel style = {font=\large,xshift=0ex},
	xticklabel style = {font=\large,yshift=0ex},
	legend style={at={(0.7,0.22)},anchor=west,font=\Large}]
\addplot [red, no marks, thick] table[x index=1, y index=2]{aster/compa_qoI_diffmesh/dict_CABH2_21.txt};
\addplot [blue, thick] table[x index=1, y index=2]{aster/compa_qoI_diffmesh/dict_CABH2_24.txt};
\addplot [orange, thick] table[x index=1, y index=2]{aster/compa_qoI_diffmesh/dict_CABH2_27.txt};
\end{axis}
\end{tikzpicture}
}
~~~
%cost EQP
\subfloat[]{
\begin{tikzpicture}[scale=0.5]
\begin{axis}[
	xlabel = {\LARGE {$t$ (s)}},
	ylabel = {\LARGE {$\varepsilon_{tt}$ - E}},
	legend entries = {$\mu^{(1)}$, $\mu^{(2)}$, $\mu^{(3)}$},
	line width=1.2pt,
	mark size=3.0pt,
	ylabel style = {font=\large,yshift=1ex},
	yticklabel style = {font=\large,xshift=0ex},
	xticklabel style = {font=\large,yshift=0ex},
	legend style={at={(0.02,0.22)},anchor=west,font=\Large}]
\addplot [red, no marks, thick] table[x index=1, y index=2]{aster/compa_qoI_diffmesh/extrados_EPMTT_21.txt};
\addplot [blue, thick] table[x index=1, y index=2]{aster/compa_qoI_diffmesh/extrados_EPMTT_24.txt};
\addplot [orange, thick] table[x index=1, y index=2]{aster/compa_qoI_diffmesh/extrados_EPMTT_27.txt};
\end{axis}
\end{tikzpicture}
}
~~~
%sampled elements
\subfloat[]{
\begin{tikzpicture}[scale=0.5]
\begin{axis}[
	xlabel = {\LARGE {$t$ (s)}},
	ylabel = {\LARGE {$\varepsilon_{zz}$ - E}},
	legend entries = {$\mu^{(1)}$, $\mu^{(2)}$, $\mu^{(3)}$},
	line width=1.2pt,
	mark size=3.0pt,
	ylabel style = {font=\large,yshift=1ex},
	yticklabel style = {font=\large,xshift=0ex},
	xticklabel style = {font=\large,yshift=0ex},
	legend style={at={(0.02,0.22)},anchor=west,font=\Large}]
\addplot [red, no marks, thick] table[x index=1, y index=2]{aster/compa_qoI_diffmesh/extrados_EPMZZ_21.txt};
\addplot [blue, thick] table[x index=1, y index=2]{aster/compa_qoI_diffmesh/extrados_EPMZZ_24.txt};
\addplot [orange, thick] table[x index=1, y index=2]{aster/compa_qoI_diffmesh/extrados_EPMZZ_27.txt};
\end{axis}
\end{tikzpicture}
}
\caption{mechanical response of a NCB under external loading.  (a) normal force on a horizontal cable, (b) tangential and (c) vertical strains on the outer wall of the standard section of the containment building,
for $\eta_{\rm dc}=5 \cdot 10^{9}$ and three values of  $\kappa$ on the coarse mesh.}
\label{fig:ss_qoIcompa_param}
\end{figure}
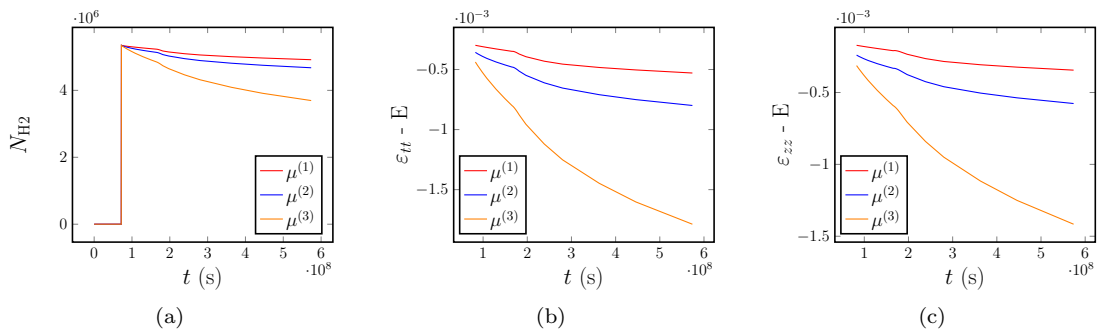

We study the solution behavior with respect to  two parameters: the desiccation creep viscosity ($\eta_{\rm dc}$) and the basic creep consolidation parameter ($\kappa$) in the parameter range
$\mu\in \mathcal{P} = [5 \cdot 10^{8}, 5\cdot 10^{10}]\times [10^{-5}, 10^{-3}] \subset \mathbb{R}^2$. 
%The state vector includes the three-dimensional displacement of the concrete structure and the one-dimensional displacement of the steel cables.
Figure \ref{fig:ss_qoIcompa_param}
shows the behavior of 
(a) the normal force on a horizontal cable, and  (b) the tangential and (c)   vertical strains on the outer wall of the standard section of the containment building, for three distinct parameter values
$\mu^{(i)} = (5.10^{9}, \kappa^{(i)})$, for $\kappa^{(i)} \in \{10^{-5}, 10^{-4}, 10^{-3}\}$, $i=1,2,3$.  Notation ``-E''  indicates that the HF data are associated  to    the outer face of the structure.
Note that  the value of the consolidation parameter $\kappa$ affects the rate of decay of the various quantities.

\subsubsection{Results}

\begin{figure}[h!]
\centering
\subfloat[]{ 
\includegraphics[width=.3\textwidth]{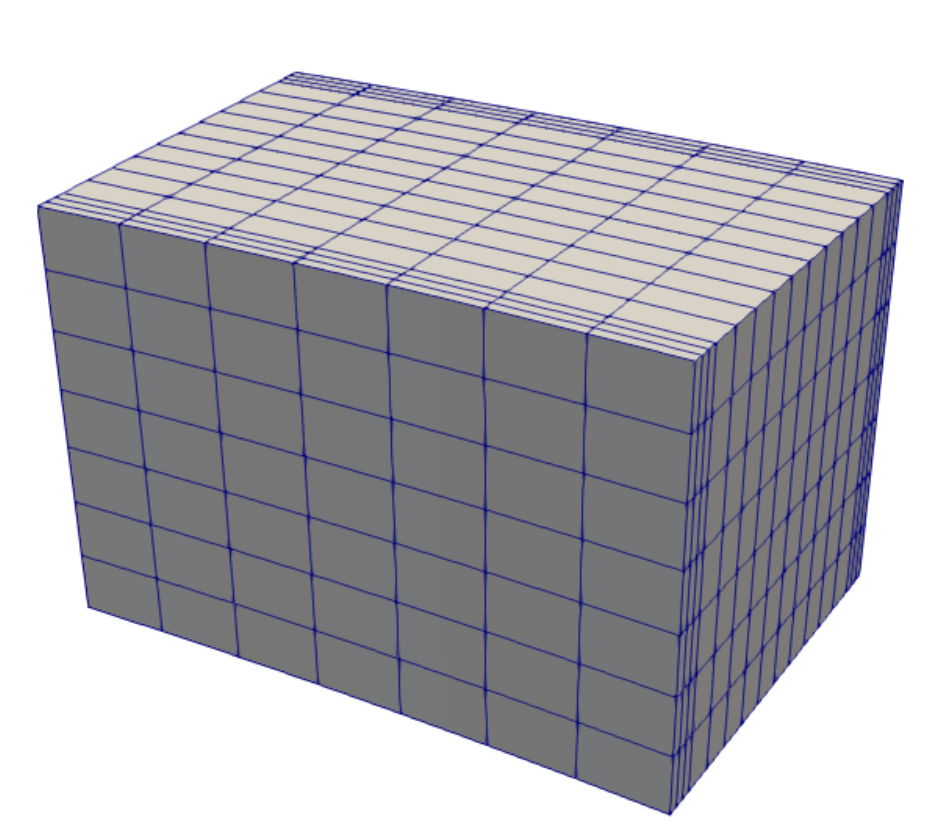}
}
\qquad
\subfloat[]{ 
\includegraphics[width=.3\textwidth]{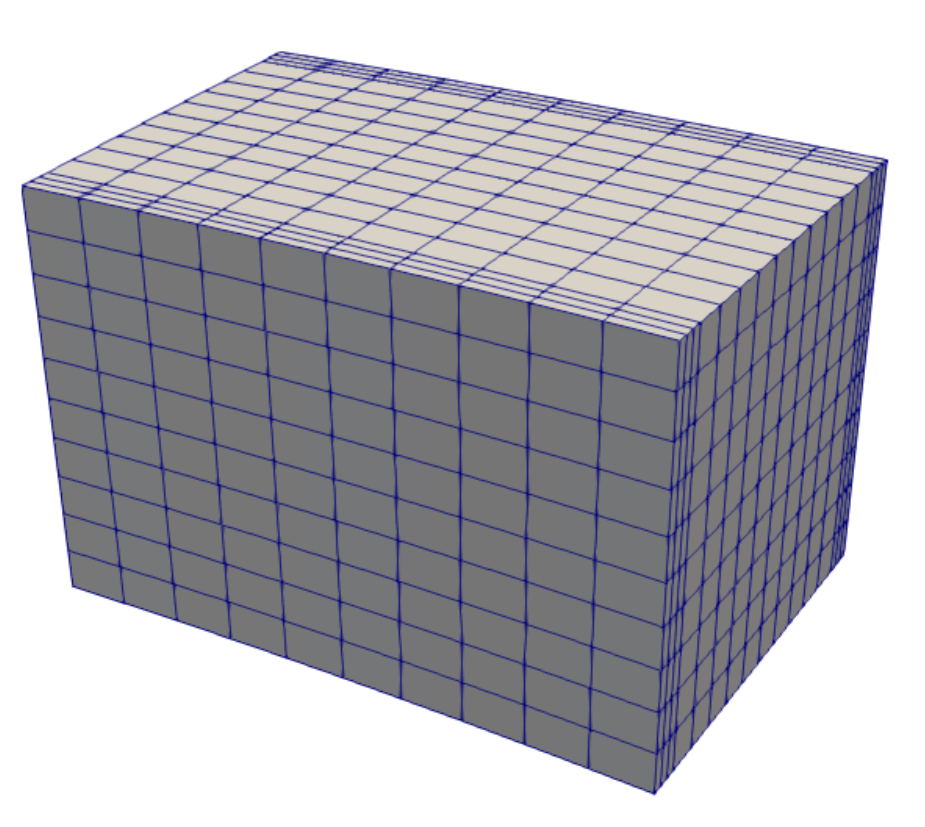}
}
\caption{mechanical response of a NCB under external loading. 
Tridimensional meshes (HEXA20) used for \texttt{code$\_$aster} calculations: (a) coarse mesh (\texttt{mesh 1}: $N_{\rm e}=784$); (b) refined mesh (\texttt{mesh 2}: $N_{\rm e}=1600$).}
\label{fig:visSS}
\end{figure}

\textbf{High-fidelity solver.} 
We consider the two distinct three-dimensional meshes, depicted in Figure \ref{fig:visSS}: a coarse mesh of $N_{\rm e}=784$ three-dimensional hexahedral elements 
and a refined mesh with $N_{\rm e}=1600$ elements. This mesh features the geometry of a portion of the building halfway up the barrel, and is crossed by two vertical and three horizontal cables. 
We consider an adaptive time-stepping scheme: approximately 45 to 50 time steps are needed  to reach the specified final time step, for all parameters in $\mathcal{P}$ and for both meshes.
Table \ref{table:ss_hf_cpu_cost} provides an overview of  the costs of the HF solver over the training set for the two meshes: we observe that the wall-clock cost of a full  HF simulation  
is roughly nine minutes for the coarse mesh and seventeen minutes for the refined mesh.
We consider a 7 by 7 training set 
$\mathcal{P}_{\rm train}$ and a 5 by 5 test set; parameters are logarithmically spaced in both directions.

\begin{table}[h!]
\centering
\begin{tabular}{|c|c|c|c|c|c|c|}
\hline
 & mean & max & min & Q1 & median & Q3 \\ \hline 
\texttt{mesh 1} & 546.91 & 905.53 & 386.96 & 387.04 & 387.12  & 387.19 \\\hline 
\texttt{mesh 2} & 1034.07 & 1658.53 & 747.60 & 748.20 & 749.78 & 749.40 \\\hline
\end{tabular}
\caption{mechanical response of a NCB under external loading. 
HF CPU cost in seconds [s] for the HF simulations on the coarse (\texttt{mesh 1}) and the refined mesh (\texttt{mesh 2})}\label{table:ss_hf_cpu_cost}
\end{table}

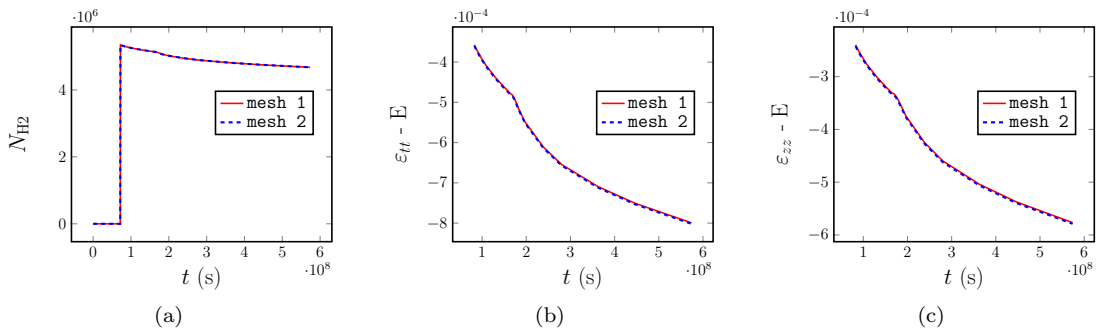
\begin{figure}[h!]
\centering
%normal effort horizontal cable
\subfloat[]{
\begin{tikzpicture}[scale=0.5]
\begin{axis}[
	xlabel = {\LARGE {$t$ (s)}},
	ylabel = {\LARGE {$N_{\rm H2}$}},
	legend entries = {\texttt{mesh 1}, \texttt{mesh 2}},
	line width=1.2pt,
	mark size=3.0pt,
	ylabel style = {font=\large,yshift=1ex},
	yticklabel style = {font=\large,xshift=0ex},
	xticklabel style = {font=\large,yshift=0ex},
	legend style={at={(0.55,0.6)},anchor=west,font=\Large}]
\addplot [red, no marks] table[x index=1, y index=2]{aster/compa_qoI_diffmesh/dict_CABH2_24.txt};
\addplot [dashed, blue, ultra thick] table[x index=1, y index=3]{aster/compa_qoI_diffmesh/dict_CABH2_24.txt};
\end{axis}
\end{tikzpicture}
}
~~~
%cost EQP
\subfloat[]{
\begin{tikzpicture}[scale=0.5]
\begin{axis}[
	xlabel = {\LARGE {$t$ (s)}},
	ylabel = {\LARGE {$\varepsilon_{tt}$ - E}},
	legend entries = {\texttt{mesh 1}, \texttt{mesh 2}},
	line width=1.2pt,
	mark size=3.0pt,
	ylabel style = {font=\large,yshift=1ex},
	yticklabel style = {font=\large,xshift=0ex},
	xticklabel style = {font=\large,yshift=0ex},
	legend style={at={(0.55,0.6)},anchor=west,font=\Large}]
\addplot [red, no marks] table[x index=1, y index=2]{aster/compa_qoI_diffmesh/extrados_EPMTT_24.txt};
\addplot [dashed, blue, ultra thick] table[x index=1, y index=3]{aster/compa_qoI_diffmesh/extrados_EPMTT_24.txt};
\end{axis}
\end{tikzpicture}
}
~~~
%sampled elements
\subfloat[]{
\begin{tikzpicture}[scale=0.5]
\begin{axis}[
	xlabel = {\LARGE {$t$ (s)}},
	ylabel = {\LARGE {$\varepsilon_{zz}$ - E}},
	legend entries = {\texttt{mesh 1}, \texttt{mesh 2}},
	line width=1.2pt,
	mark size=3.0pt,
	ylabel style = {font=\large,yshift=1ex},
	yticklabel style = {font=\large,xshift=0ex},
	xticklabel style = {font=\large,yshift=0ex},
	legend style={at={(0.55,0.6)},anchor=west,font=\Large}]
\addplot [red, no marks] table[x index=1, y index=2]{aster/compa_qoI_diffmesh/extrados_EPMZZ_24.txt};
\addplot [dashed, blue, ultra thick] table[x index=1, y index=3]{aster/compa_qoI_diffmesh/extrados_EPMZZ_24.txt};
\end{axis}
\end{tikzpicture}
}
\caption{mechanical response of a NCB under external loading. 
Comparison of the quantities of interest computed for the two meshes of the standard section: (a) normal force on a horizontal cable, (b) tangential and (c) vertical strains on the outer wall of the standard section of the containment building.}
\label{fig:ss_qoIcompa}
\end{figure}

Figure \ref{fig:ss_qoIcompa} showcases the evolution of normal forces in the central  horizontal cable ($N_{\rm H2}$), the  vertical ($\varepsilon_{zz}$) and the tangential ($\varepsilon_{tt}$) deformations on the outer surface of the geometry. 
Table \ref{table:ss_qoI_metric} 
shows the behavior of the maximum and average  relative errors 
\begin{equation}
\label{eq:QOI_error_explained}
E_{\mu}^{\rm max} = \frac{\max_{k=1,\ldots,K} | q_{\mu}^{\rm hf, (1)}(t^{(k)}) - q_{\mu}^{\rm hf, (2)}(t^{(k)}) |   }{
\max_{k=1,\ldots,K} | q_{\mu}^{\rm hf, (2)}(t^{(k)})| 
},
\quad
E_{\mu}^{\rm avg} = \frac{\sum_{k=1}^K \Delta t^{(k)}  | q_{\mu}^{\rm hf, (1)} (t^{(k)}) - q_{\mu}^{\rm hf, (2)}(t^{(k)})| \, dt  }{
\sum_{k=1}^K  | q_{\mu}^{\rm hf, (2)} (t^{(k)}) | 
}.
\end{equation}
for the three quantities of interest of Figure \ref{fig:ss_qoIcompa}.
We notice that the two meshes lead to nearly-equivalent results for this model problem.

\begin{table}[h!]
\centering
\begin{tabular}{|c|c|c|c|c|}
\hline
 &  &  &  &  \\  &$\max\limits_{\mu\in \mathcal{P}_{\rm train}} E^{\rm max}_\mu(\cdot)$&$ E^{\rm max}_\mu(\cdot)/n_{\rm train}$&$\max\limits_{\mu\in \mathcal{P}_{\rm train}} E^{\rm avg}_\mu(\cdot)$&  $ E^{\rm avg}_\mu(\cdot)/n_{\rm train}$\\ &  &  &  &  \\ \hline
 $N_{\rm H2}$ & $2.86\cdot 10^{-2}$ & $6.55\cdot 10^{-4}$ & $2.86\cdot 10^{-2}$ &  $3.91\cdot 10^{-5}$ \\  \hline
 $\varepsilon_{tt}-E$ & $4.42\cdot 10^{-2}$ & $7.68\cdot 10^{-3}$ & $7.87\cdot 10^{-3}$ & $4.84\cdot 10^{-3}$ \\  \hline
 $\varepsilon_{zz}-E$ & $5.03\cdot 10^{-2}$ &$1.09\cdot 10^{-2}$  & $1.38\cdot 10^{-2}$ & $6.73\cdot 10^{-3}$ \\  \hline
\end{tabular}
\caption{mechanical response of a NCB under external loading.  Computation of average (column 1 and 3) and maximum (column 2 and 4) errors over the training set for several errors on the quantities of interest: maximum error over all time steps (cf. \eqref{eq:QOI_error_explained}).}\label{table:ss_qoI_metric}
\end{table}

\textbf{Model reduction.}
We assess the performance of the Galerkin ROM over training and test sets. Given the sequence of parameters $\{ \mu^{\star, it} \}_{it=1}^{\texttt{maxit}}$, for $it=1,\ldots, \texttt{maxit}$, 
\begin{enumerate}
\item
we solve the HF problem to find the trajectory $\{  u_{   \mu^{\star, it}}^{(k)}  \}_{k=1}^K$;
\item
we update the reduced space $\mathcal{Z}$ using \eqref{eq:nested_space} with tolerance $tol=10^{-5}$;
\item
we update the quadrature rule using the (incremental) strategy described in section \ref{sec:unsteady_problems}.
 \end{enumerate}
Below, we assess (i) the effectiveness of the incremental strategy for the construction of the quadrature rule, and (ii) the impact of the sampling strategy on performance. Towards this end, 
we compute the projection error
\begin{equation}
\label{eq:proj_error_unsteady}
E_{it}^{\rm proj}(\mathcal{P}_{\bullet}) = 
\max_{\mu \in \mathcal{P}_{\bullet}}
\frac{\vertiii{ \mathbb{u}_{\mu}  - 
\Pi_{\mathcal{Z}_{it}}
 \mathbb{u}_{\mu}}}{\vertiii{ \mathbb{u}_{\mu}}},
\quad
\bullet \in \{\rm train, \rm test \}
\end{equation}
where $\Pi_{\mathcal{Z}_{it}}
 \mathbb{u}_{\mu} := \{  \Pi_{\mathcal{Z}_{it}} {u}_{\mu}^{(k)}  \}  $ and 
$\vertiii{ \mathbb{v} } = \sqrt{
\sum_{k=1}^K \;
\Delta t^{(k)}
\|   v^{(k)} \|^2
}$ is the discrete $L^2(0, T; \|\cdot \|)$ norm --- as in 
\cite{agouzal2024projection}, we consider $\| \bullet \| = \| \bullet \|_2$.
Further results on the prediction error and online speedups of the ROM are provided in 
\cite{agouzal2024projection}.

\begin{figure}[h!]
\centering
%nbr its
\subfloat[$\delta=10^{-4}$]{
\begin{tikzpicture}[scale=0.6]
\begin{axis}[
	xlabel = {\LARGE {$it$}},
	ylabel = {\LARGE {nbr its EQP}},
	legend entries = {incr, std},
	line width=1.2pt,
	mark size=3.0pt,
	ylabel style = {font=\large,yshift=1ex},
	yticklabel style = {font=\large,xshift=0ex},
	xticklabel style = {font=\large,yshift=0ex},
	legend style={at={(0.01,0.6)},anchor=west,font=\Large}]
\addplot [mark options={solid}, mark=square, dashed, red, mark size=2.5pt] table[x index=0, y index=6]{aster/data_ss_7x7/EQP_strat/mesh_2/nit_15_eps_1em5/delta_1em4/arr_data_track_throughiters_incrAS.txt};
\addplot [mark options={solid}, mark=triangle, dashed, blue, mark size=3.5pt] table[x index=0, y index=6]{aster/data_ss_7x7/EQP_strat/mesh_2/nit_15_eps_1em5/delta_1em4/arr_data_track_throughiters_ref.txt};
\end{axis}
\end{tikzpicture}
}
~~~
%cost EQP
\subfloat[$\delta=10^{-4}$]{
\begin{tikzpicture}[scale=0.6]
\begin{axis}[
	xlabel = {\LARGE {$it$}},
 	ylabel = {\LARGE {cost EQP [s]}},
	legend entries = {incr, std},
	line width=1.2pt,
	mark size=3.0pt,
	ylabel style = {font=\large,yshift=1ex},
	yticklabel style = {font=\large,xshift=0ex},
	xticklabel style = {font=\large,yshift=0ex},
	legend style={at={(0.01,0.6)},anchor=west,font=\Large}]
\addplot [mark options={solid}, mark=square, dashed, red, mark size=2.5pt] table[x index=0, y index=5]{aster/data_ss_7x7/EQP_strat/mesh_2/nit_15_eps_1em5/delta_1em4/arr_data_track_throughiters_incrAS.txt};
\addplot [mark options={solid}, mark=triangle, dashed, blue, mark size=3.5pt] table[x index=0, y index=5]{aster/data_ss_7x7/EQP_strat/mesh_2/nit_15_eps_1em5/delta_1em4/arr_data_track_throughiters_ref.txt};
\end{axis}
\end{tikzpicture}
}

%sampled elements
\subfloat[$\delta=10^{-4}$]{
\begin{tikzpicture}[scale=0.6]
\begin{axis}[
	xlabel = {\LARGE {$n$}},
	ylabel = {\LARGE {$\%$ sampled weights}},
	legend entries = {incr, std},
	line width=1.2pt,mark size=3.0pt,
	ylabel style = {font=\large,yshift=1ex},
	yticklabel style = {font=\large,xshift=0ex},
	xticklabel style = {font=\large,yshift=0ex},
	legend style={at={(0.01,0.6)},anchor=west,font=\Large}]
\addplot [mark options={solid}, mark=square, dashed, red, mark size=2.5pt] table[x index=0, y index=1]{aster/data_ss_7x7/EQP_strat/mesh_2/nit_15_eps_1em5/delta_1em4/arr_data_track_throughiters_incrAS.txt};
\addplot [mark options={solid}, mark=triangle, dashed, blue, mark size=3.5pt] table[x index=0, y index=1]{aster/data_ss_7x7/EQP_strat/mesh_2/nit_15_eps_1em5/delta_1em4/arr_data_track_throughiters_ref.txt}; 
\end{axis}
\end{tikzpicture}
}
~~~
\subfloat[]{
\begin{tikzpicture}[scale=0.6]
\begin{axis}[
	xlabel near ticks,
	ylabel near ticks,	
	ytick ={1,2,3,4,5,6,7,8,9,10,11},
	grid style={dashed,gray},%ymode=log,
    legend pos = north west,	
	xlabel = {\LARGE {$it$}},
 	ylabel = {\LARGE {speedups EQP}},
	line width=1.2pt,
	tick label style={font=\large},
	label style={font=\large}]
\addplot [mark options={solid}, mark=square, dashed, red, mark size=3.5pt] table[x index=0, y index=2]{aster/data_ss_7x7/EQP_strat/mesh_2/nit_15_eps_1em5/delta_1em2/speedups_incrdd.txt};
\addplot [mark options={solid}, mark=triangle, dashed, blue, mark size=3.5pt] table[x index=0, y index=2]{aster/data_ss_7x7/EQP_strat/mesh_2/nit_15_eps_1em5/delta_1em4/speedups_incrdd.txt};
\addplot [mark options={solid}, mark=otimes, dashed, orange, mark size=2.5pt] table[x index=0, y index=2]{aster/data_ss_7x7/EQP_strat/mesh_2/nit_15_eps_1em5/delta_1em6/speedups_incrdd.txt};
\legend{$\delta=10^{-2}$, $\delta=10^{-4}$, $\delta=10^{-6}$}
\end{axis}
\end{tikzpicture}
}

\caption{mechanical response of a NCB under external loading.  Progressive construction of the quadrature rule on \texttt{mesh 2}.
(a)-(b)-(c) performance of standard (``std'') and incremental (``incr'') EQ procedures for several iterations of the greedy procedures, for $\delta=10^{-4}$.
(d) computational speedup for three choices of the tolerance  $\delta$.}
\label{fig:outputs_EQP_unsteady}
\end{figure}
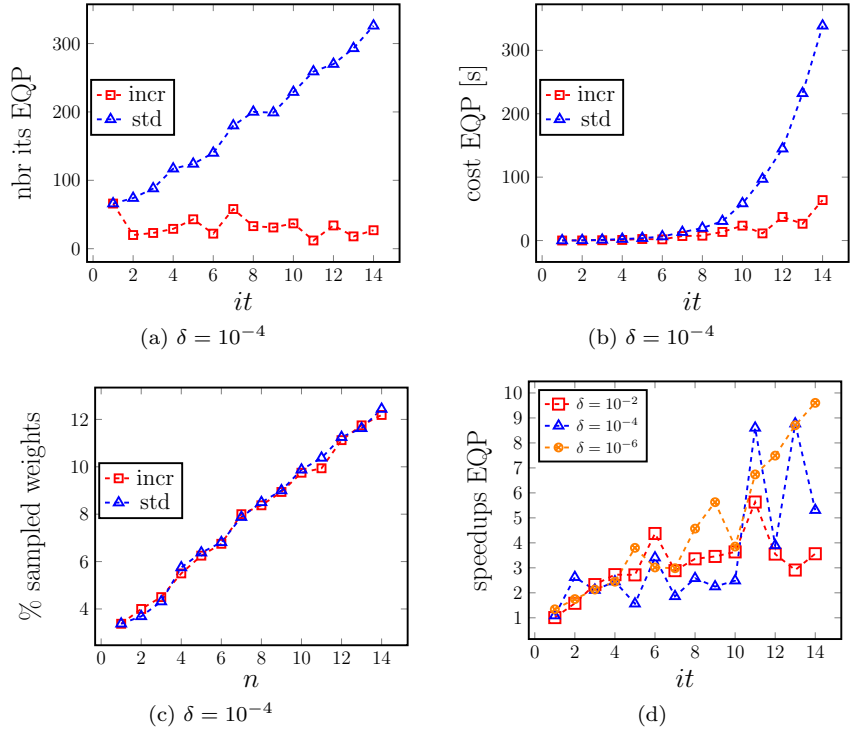

Figure \ref{fig:outputs_EQP_unsteady} illustrates the performance of the EQ procedure; in this test, we consider
the finer mesh (\texttt{mesh 2}), and we select 
 the 
parameters $\{ \mu^{\star, it} \}_{it=1}^{\texttt{maxit}=15}$ 
  using  the POD-strong greedy algorithm \ref{alg:pod_strong_greedy} based on the HF results on \texttt{mesh 1}.
Figures \ref{fig:outputs_EQP_unsteady}(a)-(b)-(c) show the number of iterations required by  Algorithm \ref{alg:nnls} to meet the convergence criterion, the computational cost, and the percentage of sampled elements, which is directly related to the online costs, for $\delta=10^{-4}$.
As for the previous test case, we observe that the progressive construction of the quadrature rule drastically reduces the NNLS iterations without hindering performance. 
Figure  \ref{fig:outputs_EQP_unsteady}(d) shows the speedup of the incremental method for three choices of the tolerance $\delta$ and for several iterations of the iterative procedure. We notice that the speedup increases as the iteration count increases:  this can be explained by observing that the percentage of new columns added during the $it$-th step in the matrix $\mathbf{G}$ decays with $it$ (cf. \eqref{eq:matrixG_unsteady}).
We further observe that the speedup increase as we decrease the tolerance $\delta$: this is justified by the fact that the number of  iterations required by Algorithm  \ref{alg:nnls} increases as we decrease 
$\delta$.

\begin{figure}[h!]
\centering
\subfloat[out-of-sample]{ 
\begin{tikzpicture}[scale = 0.5]
\begin{axis}[
	xlabel = {\LARGE {$it$}},
	ylabel = {\LARGE {proj error}},
	legend entries = {\texttt{mesh 1}, \texttt{mesh 2}, random},
	ymode=log,
	line width=1.2pt,
	mark size=3.0pt,
	ylabel style = {font=\large,yshift=1ex},
	yticklabel style = {font=\large,xshift=0ex},
	xticklabel style = {font=\large,yshift=0ex},
	legend style={at={(0.55,0.75)},anchor=west,font=\Large}]
	
\addplot [mark options={solid}, mark=square, dashed, red, mark size=2.5pt] table[x index=0, y index=2]{aster/proj_err/rel_errors/sampl_mesh1_test.txt};
\addplot [mark options={solid}, mark=triangle, dashed, blue, mark size=3.5pt] table[x index=0, y index=2]{aster/proj_err/rel_errors/sampl_mesh2_test.txt};
\addplot [mark options={solid}, dashed, black, mark=*, mark size=1.5pt] table[x index=0, y index=2]{aster/proj_err/rel_errors/sampl_ranom_1_test.txt};
\addplot [mark options={solid}, dashed, black, mark=*, mark size=1.5pt] table[x index=0, y index=2]{aster/proj_err/rel_errors/sampl_ranom_2_test.txt};
\addplot [mark options={solid}, dashed, black, mark=*, mark size=1.5pt] table[x index=0, y index=2]{aster/proj_err/rel_errors/sampl_ranom_3_test.txt};
\addplot [mark options={solid}, dashed, black, mark=*, mark size=1.5pt] table[x index=0, y index=2]{aster/proj_err/rel_errors/sampl_ranom_4_test.txt};
\addplot [mark options={solid}, dashed, black, mark=*, mark size=1.5pt] table[x index=0, y index=2]{aster/proj_err/rel_errors/sampl_ranom_5_test.txt};
\end{axis}
\end{tikzpicture}
}
\subfloat[in-sample]{ 
\centering
\begin{tikzpicture}[scale = 0.5]
\begin{axis}[
	xlabel = {\LARGE {$it$}},
	ylabel = {\LARGE {proj error}},
	legend entries = {\texttt{mesh 1}, \texttt{mesh 2}, random},
	ymode=log,
	line width=1.2pt,
	mark size=3.0pt,
	ylabel style = {font=\large,yshift=1ex},
	yticklabel style = {font=\large,xshift=0ex},
	xticklabel style = {font=\large,yshift=0ex},
	legend style={at={(0.55,0.75)},anchor=west,font=\Large}]
\addplot [mark options={solid}, mark=square, dashed, red, mark size=2.5pt] table[x index=0, y index=2]{aster/proj_err/rel_errors/sampl_mesh1_train.txt};
\addplot [mark options={solid}, mark=triangle, dashed, blue, mark size=3.5pt] table[x index=0, y index=2]{aster/proj_err/rel_errors/sampl_mesh2_train.txt};
\addplot [mark options={solid}, dashed, black, mark=*, mark size=1.5pt] table[x index=0, y index=2]{aster/proj_err/rel_errors/sampl_ranom_1_train.txt};
\addplot [mark options={solid}, dashed, black, mark=*, mark size=1.5pt] table[x index=0, y index=2]{aster/proj_err/rel_errors/sampl_ranom_2_train.txt};
\addplot [mark options={solid}, dashed, black, mark=*, mark size=1.5pt] table[x index=0, y index=2]{aster/proj_err/rel_errors/sampl_ranom_3_train.txt};
\addplot [mark options={solid}, dashed, black, mark=*, mark size=1.5pt] table[x index=0, y index=2]{aster/proj_err/rel_errors/sampl_ranom_4_train.txt};
\addplot [mark options={solid}, dashed, black, mark=*, mark size=1.5pt] table[x index=0, y index=2]{aster/proj_err/rel_errors/sampl_ranom_5_train.txt};
\end{axis}
\end{tikzpicture}
}
\subfloat[]{ 
\centering
\begin{tikzpicture}[scale = 0.5]
\begin{axis}[
	xlabel = {\LARGE {$it$}},
	ylabel = {\LARGE {$n_{it}$}},
	legend entries = {\texttt{mesh 1}, \texttt{mesh 2}, random},
	line width=1.2pt,
	mark size=3.0pt,
	ylabel style = {font=\large,yshift=1ex},
	yticklabel style = {font=\large,xshift=0ex},
	xticklabel style = {font=\large,yshift=0ex},
	legend style={at={(0.05,0.75)},anchor=west,font=\Large}]
\addplot [mark options={solid}, mark=square, dashed, red, mark size=2.5pt] table[x index=0, y index=1]{aster/proj_err/rel_errors/sampl_mesh1_train.txt};
\addplot [mark options={solid}, mark=triangle, dashed, blue, mark size=3.5pt] table[x index=0, y index=1]{aster/proj_err/rel_errors/sampl_mesh2_train.txt};
\addplot [mark options={solid}, dashed, black, mark=*, mark size=1.5pt] table[x index=0, y index=1]{aster/proj_err/rel_errors/sampl_ranom_1_train.txt};
\addplot [mark options={solid}, dashed, black, mark=*, mark size=1.5pt] table[x index=0, y index=1]{aster/proj_err/rel_errors/sampl_ranom_2_train.txt};
\addplot [mark options={solid}, dashed, black, mark=*, mark size=1.5pt] table[x index=0, y index=1]{aster/proj_err/rel_errors/sampl_ranom_3_train.txt};
\addplot [mark options={solid}, dashed, black, mark=*, mark size=1.5pt] table[x index=0, y index=1]{aster/proj_err/rel_errors/sampl_ranom_4_train.txt};
\addplot [mark options={solid}, dashed, black, mark=*, mark size=1.5pt] table[x index=0, y index=1]{aster/proj_err/rel_errors/sampl_ranom_5_train.txt};
\end{axis}
\end{tikzpicture}
}\\
\caption{mechanical response of a NCB under external loading.
Behavior of the projection error $E_{it}^{{\rm proj}}$ 
\eqref{eq:proj_error_unsteady} for  parameters selected by 
Algorithm \ref{alg:pod_strong_greedy} based on coarse (\texttt{mesh 1}) and fine data (\texttt{mesh 2});
comparison with random sampling.
 (a) performance on $\mathcal{P}_{\rm test}$ ($5 \times 5$). 
 (b) performance on $\mathcal{P}_{\rm train}$ ($7 \times 7$).
  (c)  behavior of the basis size $n_{it}$.
}
\label{fig:ss_vis_projerr}
\end{figure}
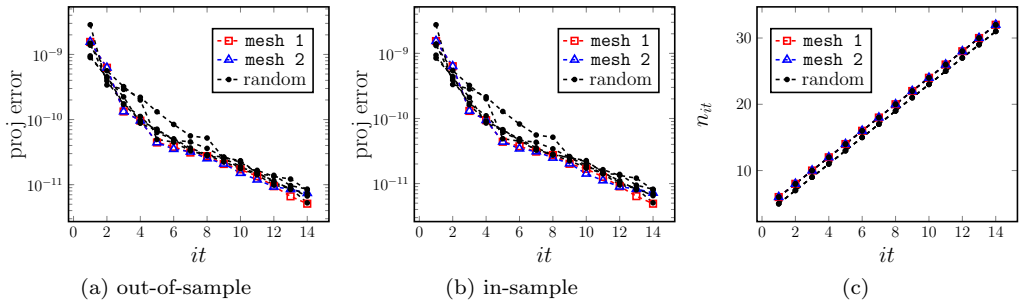

Figures \ref{fig:ss_vis_sampling} and \ref{fig:ss_vis_projerr} investigate  the efficacy of the greedy sampling strategy.
Figure \ref{fig:ss_vis_sampling} shows the parameters $\{ \mu^{\star,j} \}_j$ selected  by Algorithm \ref{alg:pod_strong_greedy} for (a) the coarse mesh (\texttt{mesh 1}) and (b) the refined mesh (\texttt{mesh 2}): we observe that the majority of the sampled parameters are 
clustered in the bottom left corner of the parameter domain for both meshes.
Figure \ref{fig:ss_vis_projerr}  shows the performance (measured through  the projection 
error \eqref{eq:proj_error_unsteady}) of the two samples depicted in Figure \ref{fig:ss_vis_sampling}  on training and test sets; to provide a concrete reference, we also compare performance with five randomly-generated samples.
Interestingly, we observe that 
for this model problem the choice of the sampling strategy has little effects on performance; nevertheless, also in this case, the greedy procedure  based on the coarse mesh consistently outperforms random sampling.

\begin{figure}[h!]
\centering
\subfloat[]{ 
\begin{tikzpicture}[scale = 0.55]
\begin{axis}[
	xlabel near ticks,
	ylabel near ticks,
	grid = both,
	xmode=log, ymode=log,
	grid style={dashed,gray},
	ylabel = {\LARGE $\kappa$},
	xlabel = {\LARGE $\eta_{\rm dc}$},
	legend pos = north east,
	tick label style={font=\large},
	label style={font=\large}]
\addplot [mark options={solid}, mark=*, gray, mark size=2pt, only marks] table[x index=1, y index=2]{aster/data_ss_7x7/sampling/mesh_1/params_txt.txt};
\addplot [mark options={solid}, mark=x, blue, mark size=3.5pt, only marks, ultra thick, nodes near coords,point meta=explicit, every node near coord/.append style={anchor=south west}] table[x index=1, y index=2, meta expr =\thisrowno{0}]{aster/data_ss_7x7/sampling/mesh_1/nit_15_eps_1em5/arr_data_indexes.txt};
\end{axis}
\end{tikzpicture}
}
\qquad
\subfloat[]{ 
\centering
\begin{tikzpicture}[scale = 0.55]
\begin{axis}[
	xlabel near ticks,
	ylabel near ticks,
	grid = both,
	xmode=log, ymode=log,
	grid style={dashed,gray},
	ylabel = {\LARGE $\kappa$},
	xlabel = {\LARGE $\eta_{\rm dc}$},
	legend pos = north east,
	tick label style={font=\large},
	label style={font=\large}]
\addplot [mark options={solid}, mark=*, gray, mark size=2pt, only marks] table[x index=1, y index=2]{aster/data_ss_7x7/sampling/mesh_2/params_txt.txt};
\addplot [mark options={solid}, mark=x, blue, mark size=3.5pt, only marks, ultra thick, nodes near coords,point meta=explicit, every node near coord/.append style={anchor=south west}] table[x index=1, y index=2, meta expr =\thisrowno{0}]{aster/data_ss_7x7/sampling/mesh_2/nit_15_eps_1em5/arr_data_indexes.txt};
\end{axis}
\end{tikzpicture}
}
\caption{mechanical response of a NCB under external loading. Parameters $\{ \mu^{\star,j} \}_j$ selected  by Algorithm \ref{alg:pod_strong_greedy} for (a) the coarse mesh (\texttt{mesh 1}) and (b) the refined mesh (\texttt{mesh 2}).
}
\label{fig:ss_vis_sampling}
\end{figure}
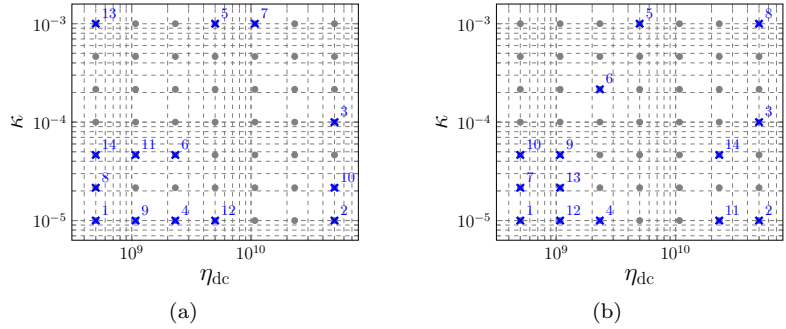

\section{Conclusions}
\label{sec:conclusions}
The computation burden of the offline training stage  remains an 
outstanding challenge of MOR techniques for nonlinear, non-parametrically-affine problems.
Adaptive sampling of the parameter domain based on greedy methods 
might contribute to reduce the number of offline HF solves that is needed to meet the target  accuracy; however, greedy methods are inherently sequential and introduce non-negligible overhead that might ultimately hinder the benefit of adaptive sampling.
To address these issues, in this work, we proposed two new strategies to accelerate 
greedy methods:  first, a progressive construction of the ROM based on HAPOD to speed up the construction of the empirical test space for LSPG ROMs and on a warm start of the NNLS algorithm to determine the empirical  quadrature rule;
second, a two-fidelity sampling strategy to reduce the number of expensive greedy iterations. 

The numerical results of section \ref{sec:numerics}
 illustrate the effectiveness and the generality of our methods
 for both steady and  unsteady problems.
 First, we found that the warm start of the NNLS algorithm enables a non-negligible reduction of the computational cost without hindering performance.
Second, we observed that the  sampling strategy based on coarse data leads to near-optimal performance: this result suggests that multi-fidelity algorithms might be particularly effective to explore the parameter domain in  the MOR training phase.

The empirical findings of this work motivate further theoretical and numerical 
investigations that we wish to pursue
in the future.
First, we wish to analyze the performance of multi-fidelity sampling methods for MOR: the ultimate goal is to devise \emph{a priori} and \emph{a posteriori} indicators to drive the choice of the mesh hierarchy and the cardinality $n_0$  of the set $\mathcal{P}_\star$ in Algorithm \ref{alg:multifidelity_greedy}.  
Second, we plan to extend the two elements of our formulation --- progressive ROM generation and multi-fidelity sampling --- to optimization problems for which the primary objective of model reduction is to estimate a suitable quantity of interest (\emph{goal-oriented MOR}): in this respect, we envision to combine our formulation with adaptive techniques for optimization
\cite{zahr2015progressive,alexandrov2001approximation,yano2021globally}.
 
\paragraph{Acknowledgements.}
This work was partly funded by ANRT (French National Association for Research and Technology) and EDF. 
The authors thank
Dr. Jean-Philippe Argaud, Dr. Guilhem 
Fert{\'e} (EDF R{\&}D) and Dr. Michel Bergmann (Inria Bordeaux) 
for fruitful discussions.
The first author thanks the code aster development team for fruitful discussions on the FE software employed for the numerical simulations of section \ref{sec:NCB}.

\paragraph{Data availability statement.}
The data that support the findings of this study are available from the corresponding author, TT, upon reasonable request.
The software developed for the investigations of section \ref{sec:NCB} is expected to be integrated in  the open-source \texttt{Python} library 
\texttt{Mordicus} \cite{mordicus2022}, which is funded by a ``French Fonds Unique
Interministériel'' (FUI) project.

\bibliographystyle{abbrv}	
\bibliography{all_refs}

\begin{thebibliography}{10}

\bibitem{aster}
{\em Finite element {code$\_$aster}, Analysis of Structures and Thermomechanics
  for Studies and Research}.
\newblock Electricit{\'e} de France (EDF), Open source on www.code-aster.org,
  1989-2024.

\bibitem{mordicus2022}
{\em Mordicus Python package. Consortium of the FUI project MOR DICUS}.
\newblock Electricit{\'e} de France (EDF), Open source on
  https://gitlab.com/mor dicus/mordicus, 2022.

\bibitem{agouzal2023projection}
E.~Agouzal, J.-P. Argaud, M.~Bergmann, G.~Fert{\'e}, and T.~Taddei.
\newblock A projection-based reduced-order model for parametric quasi-static
  nonlinear mechanics using an open-source industrial code.
\newblock {\em International Journal for Numerical Methods in Engineering,
  accepted}, 2023.

\bibitem{agouzal2024projection}
E.~Agouzal, J.-P. Argaud, M.~Bergmann, G.~Fert{\'e}, and T.~Taddei.
\newblock Projection-based model order reduction for prestressed concrete with
  an application to the standard section of a nuclear containment building.
\newblock {\em arXiv preprint arXiv:2401.05098}, 2024.

\bibitem{alexandrov2001approximation}
N.~M. Alexandrov, R.~M. Lewis, C.~R. Gumbert, L.~L. Green, and P.~A. Newman.
\newblock Approximation and model management in aerodynamic optimization with
  variable-fidelity models.
\newblock {\em Journal of Aircraft}, 38(6):1093--1101, 2001.

\bibitem{barral2023registration}
N.~Barral, T.~Taddei, and I.~Tifouti.
\newblock Registration-based model reduction of parameterized {PDE}s with
  spatio-parameter adaptivity.
\newblock {\em Journal of Computational Physics}, page 112727, 2023.

\bibitem{barrault2004empirical}
M.~Barrault, Y.~Maday, N.~C. Nguyen, and A.~T. Patera.
\newblock An ‘empirical interpolation’method: application to efficient
  reduced-basis discretization of partial differential equations.
\newblock {\em Comptes Rendus Mathematique}, 339(9):667--672, 2004.

\bibitem{benaceur2018progressive}
A.~Benaceur, V.~Ehrlacher, A.~Ern, and S.~Meunier.
\newblock A progressive reduced basis/empirical interpolation method for
  nonlinear parabolic problems.
\newblock {\em SIAM Journal on Scientific Computing}, 40(5):A2930--A2955, 2018.

\bibitem{bouhjiti2018analyse}
D.~Bouhjiti.
\newblock Analyse probabiliste de la fissuration et du confinement des grands
  ouvrages en b{\'e}ton arm{\'e} et pr{\'e}contraint-application aux enceintes
  de confinement des r{\'e}acteurs nucl{\'e}aires (cas de la maquette vercors).
\newblock {\em Academic Journal of Civil Engineering}, 36(1):464--471, 2018.

\bibitem{brand2003fast}
M.~Brand.
\newblock Fast online svd revisions for lightweight recommender systems.
\newblock In {\em Proceedings of the 2003 SIAM international conference on data
  mining}, pages 37--46. SIAM, 2003.

\bibitem{carlberg2013gnat}
K.~Carlberg, C.~Farhat, J.~Cortial, and D.~Amsallem.
\newblock The {GNAT} method for nonlinear model reduction: effective
  implementation and application to computational fluid dynamics and turbulent
  flows.
\newblock {\em Journal of Computational Physics}, 242:623--647, 2013.

\bibitem{casenave2020nonintrusive}
F.~Casenave, N.~Akkari, F.~Bordeu, C.~Rey, and D.~Ryckelynck.
\newblock A nonintrusive distributed reduced-order modeling framework for
  nonlinear structural mechanics—application to elastoviscoplastic
  computations.
\newblock {\em International journal for numerical methods in engineering},
  121(1):32--53, 2020.

\bibitem{chapman2017accelerated}
T.~Chapman, P.~Avery, P.~Collins, and C.~Farhat.
\newblock Accelerated mesh sampling for the hyper reduction of nonlinear
  computational models.
\newblock {\em International Journal for Numerical Methods in Engineering},
  109(12):1623--1654, 2017.

\bibitem{cohen2015approximation}
A.~Cohen and R.~DeVore.
\newblock Approximation of high-dimensional parametric {PDE}s.
\newblock {\em Acta Numerica}, 24:1--159, 2015.

\bibitem{conti2023multi}
P.~Conti, M.~Guo, A.~Manzoni, A.~Frangi, S.~L. Brunton, and J.~N. Kutz.
\newblock Multi-fidelity reduced-order surrogate modeling.
\newblock {\em arXiv preprint arXiv:2309.00325}, 2023.

\bibitem{du2022efficient}
E.~Du and M.~Yano.
\newblock Efficient hyperreduction of high-order discontinuous {G}alerkin
  methods: element-wise and point-wise reduced quadrature formulations.
\newblock {\em Journal of Computational Physics}, 466:111399, 2022.

\bibitem{farhat2014dimensional}
C.~Farhat, P.~Avery, T.~Chapman, and J.~Cortial.
\newblock Dimensional reduction of nonlinear finite element dynamic models with
  finite rotations and energy-based mesh sampling and weighting for
  computational efficiency.
\newblock {\em International Journal for Numerical Methods in Engineering},
  98(9):625--662, 2014.

\bibitem{farhat2015structure}
C.~Farhat, T.~Chapman, and P.~Avery.
\newblock Structure-preserving, stability, and accuracy properties of the
  energy-conserving sampling and weighting method for the hyper reduction of
  nonlinear finite element dynamic models.
\newblock {\em International journal for numerical methods in engineering},
  102(5):1077--1110, 2015.

\bibitem{feng2023multi}
L.~Feng, L.~Lombardi, G.~Antonini, and P.~Benner.
\newblock Multi-fidelity error estimation accelerates greedy model reduction of
  complex dynamical systems.
\newblock {\em International journal for numerical methods in engineering},
  124(3):5312--5333, 2023.

\bibitem{haasdonk2017reduced}
B.~Haasdonk.
\newblock Reduced basis methods for parametrized {PDEs}--a tutorial
  introduction for stationary and instationary problems.
\newblock {\em Model reduction and approximation: theory and algorithms},
  15:65, 2017.

\bibitem{haasdonk2008reduced}
B.~Haasdonk and M.~Ohlberger.
\newblock Reduced basis method for finite volume approximations of parametrized
  linear evolution equations.
\newblock {\em ESAIM: Mathematical Modelling and Numerical Analysis},
  42(2):277--302, 2008.

\bibitem{hernandez2017dimensional}
J.~A. Hernandez, M.~A. Caicedo, and A.~Ferrer.
\newblock Dimensional hyper-reduction of nonlinear finite element models via
  empirical cubature.
\newblock {\em Computer methods in applied mechanics and engineering},
  313:687--722, 2017.

\bibitem{himpe2018hierarchical}
C.~Himpe, T.~Leibner, and S.~Rave.
\newblock Hierarchical approximate proper orthogonal decomposition.
\newblock {\em SIAM Journal on Scientific Computing}, 40(5):A3267--A3292, 2018.

\bibitem{iollo2022adaptive}
A.~Iollo, G.~Sambataro, and T.~Taddei.
\newblock An adaptive projection-based model reduction method for nonlinear
  mechanics with internal variables: Application to thermo-hydro-mechanical
  systems.
\newblock {\em International Journal for Numerical Methods in Engineering},
  123(12):2894--2918, 2022.

\bibitem{lawson1995solving}
C.~L. Lawson and R.~J. Hanson.
\newblock {\em Solving least squares problems}.
\newblock SIAM, 1995.

\bibitem{MATLAB:2022}
MATLAB.
\newblock {\em R2022a}.
\newblock The MathWorks Inc., Natick, Massachusetts, 2022.

\bibitem{paul2015adaptive}
A.~Paul-Dubois-Taine and D.~Amsallem.
\newblock An adaptive and efficient greedy procedure for the optimal training
  of parametric reduced-order models.
\newblock {\em International Journal for Numerical Methods in Engineering},
  102(5):1262--1292, 2015.

\bibitem{quarteroni2015reduced}
A.~Quarteroni, A.~Manzoni, and F.~Negri.
\newblock {\em Reduced basis methods for partial differential equations: an
  introduction}, volume~92.
\newblock Springer, 2015.

\bibitem{ryckelynck2009hyper}
D.~Ryckelynck.
\newblock Hyper-reduction of mechanical models involving internal variables.
\newblock {\em International Journal for numerical methods in engineering},
  77(1):75--89, 2009.

\bibitem{sirovich1987turbulence}
L.~Sirovich.
\newblock Turbulence and the dynamics of coherent structures. {I. C}oherent
  structures.
\newblock {\em Quarterly of applied mathematics}, 45(3):561--571, 1987.

\bibitem{taddei2023compositional}
T.~Taddei.
\newblock Compositional maps for registration in complex geometries.
\newblock {\em arXiv preprint arXiv:2308.15307}, 2023.

\bibitem{taddei2021space}
T.~Taddei and L.~Zhang.
\newblock Space-time registration-based model reduction of parameterized
  one-dimensional hyperbolic {PDE}s.
\newblock {\em ESAIM: Mathematical Modelling and Numerical Analysis},
  55(1):99--130, 2021.

\bibitem{urban2014improved}
K.~Urban and A.~Patera.
\newblock An improved error bound for reduced basis approximation of linear
  parabolic problems.
\newblock {\em Mathematics of Computation}, 83(288):1599--1615, 2014.

\bibitem{veroy2003posteriori}
K.~Veroy, C.~Prud'Homme, D.~Rovas, and A.~Patera.
\newblock A posteriori error bounds for reduced-basis approximation of
  parametrized noncoercive and nonlinear elliptic partial differential
  equations.
\newblock In {\em 16th AIAA Computational Fluid Dynamics Conference}, page
  3847, 2003.

\bibitem{volkwein2011model}
S.~Volkwein.
\newblock Model reduction using proper orthogonal decomposition.
\newblock {\em Lecture Notes, Institute of Mathematics and Scientific
  Computing, University of Graz. see http://www. uni-graz.
  at/imawww/volkwein/POD. pdf}, 1025, 2011.

\bibitem{yano2019discontinuous}
M.~Yano.
\newblock Discontinuous {G}alerkin reduced basis empirical quadrature procedure
  for model reduction of parametrized nonlinear conservation laws.
\newblock {\em Advances in Computational Mathematics}, 45(5):2287--2320, 2019.

\bibitem{yano2021globally}
M.~Yano, T.~Huang, and M.~J. Zahr.
\newblock A globally convergent method to accelerate topology optimization
  using on-the-fly model reduction.
\newblock {\em Computer Methods in Applied Mechanics and Engineering},
  375:113635, 2021.

\bibitem{yano2011importance}
M.~Yano, J.~Modisette, and D.~Darmofal.
\newblock The importance of mesh adaptation for higher-order discretizations of
  aerodynamic flows.
\newblock In {\em 20th AIAA Computational Fluid Dynamics Conference}, page
  3852, 2011.

\bibitem{yano2019lp}
M.~Yano and A.~T. Patera.
\newblock An {LP} empirical quadrature procedure for reduced basis treatment of
  parametrized nonlinear {PDE}s.
\newblock {\em Computer Methods in Applied Mechanics and Engineering},
  344:1104--1123, 2019.

\bibitem{zahr2015progressive}
M.~J. Zahr and C.~Farhat.
\newblock Progressive construction of a parametric reduced-order model for
  {PDE}-constrained optimization.
\newblock {\em International Journal for Numerical Methods in Engineering},
  102(5):1111--1135, 2015.

\end{thebibliography}

\end{document}